\documentclass[11pt]{article}
\usepackage{mathtools}
\usepackage{latexsym}
\usepackage{amsmath, amsfonts,amssymb, amsthm, euscript,makeidx,mathrsfs}
\usepackage{comment}
\usepackage{tikz-cd} 
\usepackage{bbm}
\usepackage{enumerate}
\usepackage[utf8]{inputenc}

\numberwithin{equation}{section}
\newtheorem{definition}{Definition}[section]
\newtheorem{theorem}{Theorem}[section]

\newtheorem{proposition}[theorem]{Proposition}

\newtheorem{assumption}{Assumption}[section]

\def\qed{ \ \vrule width.2cm height.2cm depth0cm\smallskip}

\newcommand{\hP}{\hat\dbP}

\newcommand{\ba}{\begin{array}}
\newcommand{\ea}{\end{array}}
\newcommand{\be}{\begin{equation}}
\newcommand{\ee}{\end{equation}}
\newcommand{\bea}{\begin{eqnarray}}
\newcommand{\eea}{\end{eqnarray}}
\newcommand{\beaa}{\begin{eqnarray*}}
\newcommand{\eeaa}{\end{eqnarray*}}

\def\neg{\negthinspace}

\def\dbE{\mathbb{E}}
\def\dbF{\mathbb{F}}

\def\dbR{\mathbb{R}}

\def\dbx{\mathbbm{x}}

\def\a{\alpha}
\def\b{\beta}

\def\d{\delta}
\def\e{\varepsilon}

\def\l{\lambda}
\def\m{\mu}
\def\n{\nu}
\def\si{\sigma}
\def\t{\tau}
\def\f{\varphi}
\def\th{\theta}
\def\o{\omega}
\def\h{\widehat}

\def\G{\Gamma}
\def\D{\Delta}
\def\Th{\Theta}
\def\L{\Lambda}

\def\O{\Omega}
\def\cA{{\cal A}}
\def\cB{{\cal B}}

\def\cF{{\cal F}}
\def\cG{{\cal G}}

\def\cN{{\cal N}}

\def\cS{{\cal S}}

\def\hC{\mathbb{C}}
\def\hD{\mathbb{D}}
\def\hE{\mathbb{E}}
\def\hF{\mathbb{F}}

\def\hN{\mathbb{N}}

\def\hP{\mathbb{P}}

\def\hR{\mathbb{R}}

\def\sA{\mathscr{A}}
\def\sB{\mathscr{B}}

\def\sD{\mathscr{D}}

\def\sJ{\mathscr{J}}

\def\sP{\mathscr{P}}

\def\sT{\mathscr{T}}
\def\sU{\mathscr{U}}

\def\no{\noindent}

\def\ss{\smallskip}
\def\ms{\medskip}
\def\bs{\bigskip}
\def\q{\quad}
\def\qq{\qquad}

\def\pa{\partial}
\def\cd{\cdot}
\def\cds{\cdots}
\def\lan{\langle}
\def\ran{\rangle}

\def\qed{ \hfill \vrule width.25cm height.25cm depth0cm\smallskip}
\newcommand{\dfnn}{\stackrel{\triangle}{=}}
\newcommand{\basa}{\begin{assumption}}
\newcommand{\easa}{\end{assumption}}

\newcommand{\ol}{\overline}
\newcommand{\ul}{\underline}

\newcommand{\bas}{\begin{assum}}
\newcommand{\eas}{\end{assum}}

\def\lan{\mathop{\langle}}
\def\ran{\mathop{\rangle}}

\def\pa{\partial}
\def\h{\widehat}
\def\wt{\widetilde}
 \def\cd{\cdot}
\def\cds{\cdots}

\def\dis{\displaystyle}
\def\wt{\widetilde}

\def\bmx{\mathbbm{x}}

\def\1{{\bf 1}}

\def\:{\!:\!}

at 9pt

\begin{document}

\newtheorem{thm}{Theorem}[section]
\newtheorem{lem}[thm]{Lemma}
\newtheorem{cor}[thm]{Corollary}
\newtheorem{prop}[thm]{Proposition}
\newtheorem{rem}[thm]{Remark}
\newtheorem{eg}[thm]{Example}
\newtheorem{defn}[thm]{Definition}
\newtheorem{assum}[thm]{Assumption}

\renewcommand {\theequation}{\arabic{section}.\arabic{equation}}
\def\thesection{\arabic{section}}

\title{\bf Equilibrium Model of Limit Order Books -- A Mean-field Game View}

\author{
Jin Ma\thanks{ \noindent Department of
Mathematics, University of Southern California, Los Angeles, CA 90089, USA.
Email: jinma@usc.edu. This author is supported in part
by US NSF grants \#DMS-1908665.}
 ~ and ~{Eunjung Noh}\thanks{\noindent
Department of Mathematics, Rutgers University, Piscataway, NJ 08854, USA. 
E-mail: en221@math.rutgers.edu 
 }}

\date{\today}


\maketitle
\begin{abstract}

In this paper we study a continuous time equilibrium model of limit order book (LOB) in which the liquidity dynamics follows a
non-local, reflected mean-field stochastic differential equation (SDE) with evolving intensity. Generalizing the basic idea of \cite{MWZ}, 
we argue that the frontier of the LOB (e.g., the best asking price) is the value function of a mean-field  stochastic control problem, 
as the limiting version of a Bertrand-type competition among the liquidity providers. With a detailed analysis on the $N$-seller static Bertrand
game, we formulate a continuous time limiting mean-field control problem of the representative seller. We then validate the dynamic programming principle (DPP), and show that the value function is a viscosity solution of the corresponding Hamilton-Jacobi-Bellman (HJB) equation. We argue
that the value function can be used to obtain the equilibrium density function of the LOB,  following the idea of \cite{MWZ}.
\end{abstract}

\vfill

\no{\bf Keywords.} \rm Limit order book, equilibrium model, Bertrand games, reflected mean-field SDEs with jumps, Hamilton-Jacobi-Bellman equation, viscosity solution.

\bs

\no{\it 2000 AMS Mathematics subject classification:} 
60H07, 60H30,
91B54, 
93E20. 

\eject


\section{Introduction}
\setcounter{equation}{0}

With the rapid growth of electronic trading, the study of order-driven markets has become an increasingly prominent focus in quantitative finance. Indeed, in the current financial world more than half of the markets use a limit order book (LOB) mechanism to facilitate trade. There has been a 
large amount of literature studying LOB from various angles, combined with some associated optimization problems such as placement, 
liquidation, executions,  etc. (see, e.g. \cite{Alfonsi1}, \cite{Alfonsi3}, \cite{Alfonsi4}, \cite{Avellaneda}, \cite{controlled intensity}, \cite{Cont}, \cite{Gatheral}, \cite{Charles}, \cite{Lokka}, \cite{Obizhaeva}, \cite{Potters}, \cite{Predoiu}, to mention a few).
Among many important structural issues of LOB, one of the focuses has been the dynamic movement of the LOB, both its frontier and its ``density" (or ``shape"). The latter was shown to be a determining factor of the ``liquidity cost" (cf. \cite{MWZ}), an important aspect that impacts
the pricing of the asset. We refer to, e.g., \cite{Alfonsi2,Nadtochiy1,4L,MWZ} for the study of  LOB particularly concerning its shape formation. 

%
%

In this paper, we try to extend dynamic model of LOB proposed in  \cite{MWZ} in two major aspects. The guiding idea is to specify the 
{\it expected equilibrium utility function}, which plays an essential role in the modeling of the shape of the LOB in that it endogenously 
determines both the dynamic density of the LOB and its frontier. More precisely, instead of assuming, more or less in an ad hoc manner,
that the equilibrium price behaves like an ``utility function", we shall consider it as the consequence of a Bertrand-type game among a 
large number of liquidity providers (sellers who set limit orders). Following the argument of \cite{CS15}, we first study an $N$-seller static
Bertrand game, each with a profit function involving not only the limit order price less the waiting cost, the same criterion as that in \cite{MWZ}, 
but also the average of the other sellers limit orders observed. We show
that the Nash equilibrium exists in such a game. With an easy randomization argument, we can then show that, as $N\to \infty$, the Nash
equilibrium converges to an optimal strategy of a single player's optimization problem with a mean-field nature, as expected.

We note that the Bertrand game in finance can be traced back to as early as 1800s, when Cournot \cite{Cournot} and Bertrand \cite{Bertrand} first studied oligopoly models of markets with a small number of competitive players. We refer to Friedman \cite{Friedman} and Vives \cite{Vives} for background and references. 
Since Cournot's model uses quantities as a strategic variable to 
determine the price, while Bertrand model does the opposition, we choose to use the Bertrand game as it fits our problem better. 
We shall assume that the sellers use the same marginal profit function, but with different choices of the price-waiting cost preference to achieve 
the optimal outcome (see \S3 for more detailed formulation). We would like to point out that our study of Bertrand game is in a sense ``motivational"
for the second main feature of this paper, that is, the continuous time, mean-field type dynamic liquidity model. More precisely, we assume that 
the liquidity dynamics is a pure-jump Markov process, with a mean-field type state dependent jump intensity. Such a dynamic game is rather complicated, and is expected to involve systems of  nonlinear, mean-field type partial differential equations (see, e.g., \cite{HHS, LS11}). We therefore consider the limiting case as the number of sellers tends to infinity,  
and argue that the dynamics of 
the  total liquidity should follow a pure jump Markov process with a mean-field type intensity, and can  be expressed
as the solution of a pure-jump SDE with reflecting boundary conditions and mean-field type state-dependent
jump intensity. We note that such SDE is itself new and therefore interesting in its own right.  
 
We should point out that the special features of our underlying liquidity dynamics 
(mean-field type; state-dependent intensity; and reflecting boundary conditions) require the combined technical tools in mean-field games, McKean-Vlasov SDEs with state-dependent jump intensities, and SDEs 
with discontinuous paths and reflecting boundary conditions. In particular, we refer to the works  \cite{MFG note,Carmona1,Carmona2,MFG application,MFjump,MFG,Ma93,MYZ} (and the references cited therein) for the technical foundation of this paper. Furthermore, 
apart from justifying the underlying liquidity dynamics, another main task of this paper is to substantiate the corresponding stochastic control 
problem, including validating the dynamic programming principle (DPP) and showing that the value function is a viscosity solution to the corresponding
Hamilton-Jacobi-Bellman (HJB) equation.

This paper is organized as follows. In \S2 we introduce necessary notations and preliminary concepts, and study the well-posedness of 
a reflected mean-field SDEs with jumps that will be essential in our study. We shall also provide 
an It\^o's formula involving reflected mean-field SDEs with jumps for ready reference. 
In \S3 we investigate a static Bertrand game with $N$ sellers, and its limiting behavior as $N$ tends to infinity. 
Based on the results, we then propose in \S4 a continuous time mean-field type stochastic control problem for a representative seller, 
as the limiting version of dynamic Bertrand game when the number of sellers becomes sufficiently large. 
In \S5 and \S6 we prove the dynamic programming principle (DPP), derive the HJB equation, and show that the value function is a viscosity solution to the corresponding HJB equation. 
Finally, in \S7 we make some concluding remarks regarding how the value function is related to the equilibrium density of the LOB, following the arguments in \cite{MWZ}.
\section{Preliminaries}

Throughout this paper we let  $(\Omega, \mathcal{F}, \mathbb{P})$ be a complete probability space on which is defined two standard Brownian
motions $W = \lbrace W_{t}: t \geq 0 \rbrace$ and $B = \lbrace B_{t}: t \geq 0 \rbrace$. Let $(\cA, \sB_{\cA})$  and $(\cB, \sB_{\cB})$ be two measurable spaces. We assume that there are 
two Poisson random measures $\mathcal{N}^{s}$ and $\cN ^b$, defined on $\hR_{+} \times \cA\times \hR_{+} $ and $\mathbb{R_{+}} \times \cB$,  and with L\'evy measures $\nu^s(dz)$ and $\nu^b(dz)$, respectively.  In other words, we assume that the Poisson measures
$\cN^s$ and $\cN^b$ have mean measures $\h\cN^s(\cd):=m\times \n^s\times m(\cd)$ and $\h\cN^b(\cd):=m\times \n^b(\cd)$, respectively, where $m(\cd)$ denotes the Lebesgue measure on $\hR_+$, and we denote the \textit{compensated} random measures
 $\wt\cN^s(A) :=(\cN^s-\h\cN^s)(A)=\cN^s(A)- (m\times \n^s\times m)(A)$ and $\wt\cN^b(B):=(\cN^b-\h\cN^b)(B)=\cN^b(B)-(m\times \n^b)(B)$, for any 
$A\in \sB(\hR_+\times \cA\times\hR_+)$ and $B\in \sB(\hR_+\times \cB)$.
For simplicity, throughout this paper we 
assume that both $\n^s$ and $\n^b$ are finite, that is, $\n^s(\cA), \n^b(\cB)< \infty$, and 
we assume the Brownian motions and Poisson random measures are mutually independent.
We note that for any $A\in\sB(\cA\times\hR_+)$ and $B\in \sB(\cB)$, the processes $(t, \o)\mapsto \wt \cN^s([0,t]\times A, \o)$, $\wt \cN^b([0,t]\times B, \o)$ are both $\hF^{\cN^s, \cN^b}$-martingales. Here $\mathbb{F}^{\cN^s, \cN^b}$ denotes
the  filtration generated by $\cN^s$ and $\cN^b$. 


For a generic Euclidean space $E$ and for $T>0$, we denote $C([0,T];E)$ and $\hD([0,T];E)$ to be the spaces of continuous
and  c\`adl\`ag functions, respectively. We endow both spaces with ``sup-norms", so that both of them are complete metric spaces. 
Next, for $p\ge 1$ we denote $L^{p}(\mathcal{F}; E)$ to be the space of all $E$-valued $\mathcal{F}$-measurable random variable $\xi$ defined
on the probability space $(\O, \cF, \hP)$ such that $\hE[|\xi|^p]<\infty$.  In particular,   $L^{2}(\mathcal{F}; \mathbb{R})$ is a Hilbert space with inner product $(\xi, \eta)_{2} = \mathbb{E}[\xi \eta], 
\; \xi, \eta \in L^{2}(\mathcal{F}; \mathbb{R})$, and a norm $\| \xi \|_{2} = (\xi, \xi)_{2}^{1/2}$. 
Also, we denote $L^{p}_{\dbF} ([t,T]; E)$ to be all $E$-valued $\hF$-adapted process $\eta$ on $[t,T]$, such that
$\| \eta \|_{p,T} := \dbE [ \int_t^T | \eta_s |^{p} ds ] ^{1/p} < \infty$. We often use the notations $L^p(\hF; C([0,T];E))$
 and $L^p(\hF; \hD([0,T];E))$
when we need to specify the path properties for elements in $L^p_\hF([0,T];E)$.

For $p\ge1$ we denote by $\sP_{p}(E)$ the space of probability measures $\mu$ on $(E, \sB(E))$ with finite $p$-th moment, i.e. $\| \mu \| _{p}^{p} := \int_{E} |x|^p \mu(dx) < \infty$. Clearly, for  $\xi\in L^p(\cF; E)$, the law $\mathcal{L}(\xi)=\mathbb{P}_{\xi} := \mathbb{P} \circ \xi^{-1}\in \sP_p(E)$. We  endow  $\sP_{p}(E)$  with the following $p$-Wasserstein metric:  
\bea
\label{W2}
W_{p}(\mu,\nu) & :=& \inf \Big\{\Big( \int_{E \times E} |x-y|^{p} \pi(dx,dy) \Big)^{\frac{1}{p}} : \pi \in \sP_{p}(E \times E) \text{\;with marginals\;} \mu \text{\;and\;} \nu \Big\} \nonumber\\
& =& \inf \Big\{ \|\xi - \xi'\|_{L^{p}( \O)}: \xi, \xi' \in L^{p}(\mathcal{F}; E) \text{\;with\;} \mathbb{P}_{\xi} = \mu, ~\mathbb{P}_{\xi'} = \nu \Big\}.
\eea

Furthermore, we suppose that there is a sub-$\sigma$-algebra $\mathcal{G} \subset \mathcal{F}$ such that 
(i)  the Brownian motion $W$ and Poisson random measures $\cN ^s, \cN ^b$ are independent of $\mathcal{G}$; and (ii)
$\mathcal{G}$ is ``rich enough" in the sense that for every 
$\mu \in \sP_{2}(\dbR)$,
 there is a random variable $\xi \in L^{2}(\mathcal{G}; E)$ such that $\mu = \mathbb{P}_{\xi}$.
Let $\mathbb{F} = \hF^{W, B,\cN^s, \cN^b \vee \cG}=\{\cF_t\}_{t\ge0}$, where $\cF_{t}=\cF^W_t\vee \cF^B_t\vee \cF_t^{\cN^s}\vee \cF^{\cN^b}_t\vee \cG$, $t \geq 0$,  be the filtration generated by $W$, $B$, $\cN^s$, $\cN^b$, and $\cG$, augmented by all the $\hP$-null sets so that it satisfies the {\it usual
hypotheses} (cf. \cite{Protter}).

The following notation of ``differentiability" with respect to probability measures is based on the lecture notes \cite{MFG note} following the course at Coll\'ege de France by P. L. Lions.
For a function $f: \sP_{2}(\mathbb{R}) \rightarrow \mathbb{R}$, we introduce a ``lift" function $f^{\sharp}: L^{2}(\mathcal{F}; \mathbb{R}) \rightarrow \mathbb{R}$ such that $f^{\sharp}(\xi) := f(\hP_{\xi}), \; \xi \in L^{2}(\mathcal{F}; \mathbb{R})$. 
Clearly $f^\sharp$ depends only on the law of $\xi \in L^{2}(\mathcal{F}; \mathbb{R})$, and is independent of the choice of the representative $\xi$. 
We say that $f: \sP_{2}(\mathbb{R}) \rightarrow \mathbb{R}$ is \textit{differentiable} at $\mu_0 \in \sP_{2}(\mathbb{R})$ if there exists $\xi_{0} \in L^{2}(\mathcal{F}; \mathbb{R})$ with $\hP_{\xi_{0}} = \mu_0$ such that $f^{\sharp}$ is Fr\'echet differentiable at $\xi_0$. In other words, there exists a continuous linear functional $Df^{\sharp}(\xi_0): L^{2} (\mathcal{F}; \mathbb{R}) \rightarrow \mathbb{R}$ such that
\begin{equation}
\label{FrechetDiff}
f^{\sharp}(\xi_{0} + \eta) - f^{\sharp}(\xi_0) = Df^{\sharp}(\xi_{0})(\eta) + o(\| \eta \|_{2}).
\end{equation} 
We shall denote $D_{\eta}f(\mu_0)=Df^{\sharp}(\xi_{0})(\eta)$, and refer to it as the Fr\'echet derivative of $f$ at $\mu_0$ in the direction $\eta$. 
By Riesz' Representation Theorem, there exists a unique random variable $\zeta \in L^{2}(\mathcal{F}; \mathbb{R})$ such that $D_{\eta}f(\mu_0) = Df^{\sharp}(\xi_{0})(\eta) = (\zeta, \eta)_{2} = \mathbb{E}[\zeta\eta]$, $\eta \in L^{2}(\mathcal{F}; \mathbb{R})$. It was shown in \cite[Lemma 3.2]{MFG note} that there exists a Borel function $h[\mu_0]: \mathbb{R} \rightarrow \mathbb{R}$, which depends only on the law $\mu_0 = \hP_{\xi_0}$  but not on the particular choice of the representative $\xi_0$, such that $\zeta = h[\mu_0](\xi_{0})$, $\hP$-a.s.. We shall denote $\partial_{\mu} f(\hP_{\xi_0}, y) \triangleq h[\mu_0](y), y \in \mathbb{R}$ and refer to it the derivative of $f: \sP_{2}(\mathbb{R}) \rightarrow \mathbb{R}$ at $\mu_0 = \hP_{\xi_0}$. In other words, we have the following identities:
\begin{equation}
Df^{\sharp}(\xi_0) = \zeta = h[\hP_{\xi_0}](\xi_0) = \partial_{\mu}f(\hP_{\xi_0}, \xi_0),
\end{equation}
and (\ref{FrechetDiff}) can be rewritten as 
\begin{equation}
f(\hP_{\xi_{0} + \eta}) - f(\hP_{\xi_0}) = \mathbb{E}[h[\hP_{\xi_{0}}](\xi_{0}) \eta] + o(\| \eta \|_{2})=(\partial_{\mu}f(\hP_{\xi_0}, \xi_0), \eta)_2+o(\|\eta\|_2),
\end{equation}
 and $ D_{\eta}f(\hP_{\xi_0}) = (\partial_{\mu}f(\hP_{\xi_0}, \xi_0), \eta)_{2}$ where $\eta = \xi - \xi_0$. Notice that $\partial_{\mu}f(\hP_{\xi_0}, y)$ is only $\hP_{\xi_0}(dy)$-a.e. uniquely determined.

Let us introduce two spaces that are useful for our analysis later. We denote $C_{b}^{1,1}(\sP_{2}(\mathbb{R}))$ the space of 
all differentiable functions $f: \sP_{2}(\mathbb{R}) \rightarrow \mathbb{R}$ such that 
$\partial_{\mu} f$ exists, and 
is {\it bounded} and {\it Lipschitz continuous}. That is, for some constant $C >0$, it holds

\ss
(i) $| \partial_{\mu} f(\mu, x) | \leq C, \;\; \mu \in \sP_{2}(\mathbb{R}), \; x \in \mathbb{R}$; 

\ss
(ii)  $| \partial_{\mu} f(\mu, x) - \partial_{\mu} f(\mu', x')| 
\leq C \lbrace |x-x'| + W_{2}(\mu,\mu') \rbrace, \;\; \mu,\mu' \in \sP_{2}(\mathbb{R}), \; x,x' \in \mathbb{R}$.  

\ms

Note that if $f \in C_{b}^{1,1}(\sP_{2}(\mathbb{R}))$, then for fixed $y \in \mathbb{R}$, we can discuss the differentiability of the 
derivative function $ \partial_{\mu} f (\cdot, y) : \sP_{2}(\mathbb{R}) \rightarrow \mathbb{R} $. In particular, if $ \partial_{\mu} f (\cdot, y) \in C_{b}^{1,1}(\sP_{2}(\mathbb{R}))$ again, then for every $y \in \mathbb{R}$, we can define
\bea
\label{f2mu}
\partial_{\mu}^{2} f(\mu, x,y) := \partial_{\mu}((\partial_{\mu}f)(\cdot, y))(\mu, x), \;\; (\mu, x,y) \in \sP_{2}(\mathbb{R}) \times \mathbb{R} \times \mathbb{R}.  
\eea
We shall  denote $C_{b}^{2,1}(\sP_{2}(\mathbb{R}))$ to be the space of all functions $f \in C_{b}^{1,1}(\sP_{2}(\mathbb{R}))$ such that 

\ss 
(i) $\partial_{\mu} f (\cdot, x) \in C_{b}^{1,1}(\sP_{2}(\mathbb{R}))$ for all $x\in \mathbb{R}$;

\ss
(ii) $\partial_{\mu}^{2} f: \sP_{2}(\mathbb{R}) \times \mathbb{R} \times \mathbb{R} \rightarrow \mathbb{R} \otimes \mathbb{R} $ is bounded and Lipschitz continuous;

\ss
(iii) $\partial_{\mu}f(\mu, \cdot) : \mathbb{R} \rightarrow \mathbb{R}$ is differentiable for every $\mu \in \sP_{2}(\mathbb{R})$, 
and its derivative $\partial_{y}\partial_{\mu} f: \sP_{2}(\mathbb{R}) \times \mathbb{R} \rightarrow \mathbb{R} \otimes \mathbb{R} $ is bounded and Lipschitz continuous. 


\subsection{Mean-field SDEs with reflecting boundary conditions}

In this subsection we consider the following (discontinuous) SDE with reflection: for $t\in[0,T)$,
\bea
\label{MFSDER}
X_{s} &=& x + \int_t^s \int_{A \times \hR_{+}}  \th(r, X_{r-},  \hP_{X_{r}},z) 
\1_{ [0, \l(r, X_{r-}, \hP_{X_{r}})] } (y) \wt\cN^{s}(dr dz dy)\\
&& +\int_t^s b(r, X_r, \hP_{X_r})dr+ \int_t^s \si(r, X_r, \hP_{X_r})dB_r + \b_s + K_s, \qq s\in[t, T], \nonumber
\eea
where $\th$, $\l$, $b$, $\sigma$ are measurable functions defined on appropriate subspaces of $[0,T]\times \O\times \hR\times \sP_2(\hR)\times \hR$, $\b$ is an $\hF$-adapted process with c\`adl\`ag paths, and $K$ is a ``reflecting process", that is, it is an
$\hF$-adapted, non-decreasing, c\`adl\`ag process, so that 

\ss
(i) $X_s\ge 0$, $\hP$-a.s.;

\ss
(ii) $\int_0^T\1_{\{X_r>0\}}dK^c_r=0$, $\hP$-a.s. ($K^c$ denotes the continuous part of $K$); and 

\ss
(iii) $\D K_t =  (X_{t-} + \D Y_t )^{-}$, where $Y=X-K$.

\ms
We call SDE (\ref{MFSDER}) a {\it mean-field SDE with discontinuous paths and reflections} (MFSDEDR), and we denote the solution
by $(X^{t,x}, K^{t,x})$, although the superscript is often omitted when context is clear. If $b,\sigma =0$ and $\b$ is pure jump, 
then the solution $(X, K)$ becomes pure jump as well (i.e., $dK^c\equiv 0$). We note that the main feature of this SDE is
that the jump intensity $\l(\cds)$ of the solution $X$ is ``state-dependent" with mean-field nature. Its well-posedness thus requires some attention since, to the best of our knowledge, it has not been studied in the literature. 
%

We shall make use of the following {\it Standing Assumptions}. 
\begin{assum}
\label{assump1}
The mappings $\lambda: [0,T]\times \mathbb{R} \times \sP_{2}(\mathbb{R}) \mapsto \mathbb{R}_+$, $b: [0,T]\times \O\times \hR \times \sP_{2}(\mathbb{R})\mapsto \hR$, $\si: [0,T]\times \O\times \hR \times \sP_{2}(\mathbb{R})\mapsto \hR$, and $\th:[0,T]\times \O\times \hR\times \sP_{2}(\mathbb{R})\times \hR\mapsto \hR$ are all uniformly bounded and continuous in $(t,x)$, 
and satisfy the following conditions, respectively:

\ss
(i) For fixed $\m\in\sP_2(\hR)$ and $x, z\in\hR$,  the mappings $(t,\o)\mapsto \th(t,\o, x, \m, z), (b, \si)(t,\o, x, \m)$ are $\hF$-predictable; 

\ss
(ii) For fixed $\m\in\sP_2(\hR)$,  $(t, z)\in[0,T]\times \hR$, and $\hP$-a.e $\o\in\O$,  the functions $ \l(t,\cd, \m)$, $b(t,\o,\cd, \m)$,
 $\si(t,\o,\cd, \m)$, $\th(t, \o, \cd, \m, z)\in C^1_b(\hR)$;

\ss
(iii) For fixed $(t, x, z)\in[0,T]\times \hR^2$, and $\hP$-a.e $\o\in\O$, 
the functions $ \l(t,x, \cd)$, $b(t,\o, x,\cd)$, $\si(t, \o, x, \cd)$, $\th(t, \o,x,  \cd, z)\in C_{b}^{1,1}(\sP_{2}(\hR))$;


\ss 
(iv) There exists $L>0$, such that  for $\hP$-a.e. $\o\in\O$, it holds that 
\beaa
&& |{\lambda}(t, x,\mu) - {\lambda}(t, x',\mu')|+ | b(t, \o, x,\m) - b(t,\o, x', \m') | \\
&& + | \si(t, \o, x,\m) - \si(t,\o, x', \m') |+| \th(t, \o, x,\m, z) - \th(t, \o, x', \m',z) |\\
&& \qq\qq \le L \left(|x-x'| + W_1(\m, \m')\right), \qq t\in[0,T], ~ x,x', z \in \mathbb{R}, ~\mu, \mu' \in \sP_{2}(\mathbb{R}).
\qq\qq \qq\qed
 \eeaa
\end{assum}

\begin{rem}
\label{Remark1_S2}{\rm
(i) The requirements on the coefficients in Assumption \ref{assump1} (such as boundedness) are stronger than necessary, only to simplify the arguments. More general (but standard) assumptions  using, e.g.,  the $L^2$-integrability with respect to the L\'evy measure $\n$ (even in the case when
$\n(\hR)=\infty$) are easily extendable without substantial difficulties. We prefer not to pursue such generality since this is not the main
purpose of the paper. 

(ii) Throughout this paper, unless specified, we shall denote $C>0$ to be a generic constant depending only on $T $, $\n(\hR) $, 
and bounds 
in Assumption \ref{assump1}. Furthermore, we shall allow it  to vary from line to line. 
\qed}
\end{rem}

It is well-known that (see, e.g., \cite{BLPR}), as a mean-field SDE, the solution to (\ref{MFSDER}) may not satisfy the so-called
``flow property", in the sense that 
$X_{r}^{t,x} \neq X_{r}^{s,X_{s}^{t,x}}, 0 \leq t \leq s \leq r \leq T$.  It is also noted in \cite{BLPR} that if we consider the  following 
accompanying SDE of (\ref{MFSDER}):
\bea
\label{Xtxi}
X_{s}^{t,\xi} &=& \xi + \int_t^s \int_{A \times \hR_{+}}  \th(r, X^{t,\xi}_{r-},  \hP_{X^{t,\xi}_{r}}, z) 
\1_{ [0, \l(r, X^{t,\xi}_{r-}, \hP_{X^{t,\xi}_{r}})] } (y) \wt\cN^{s}(dr dz dy)\\
&& +\int_t^s b(r, X^{t,\xi}_r, \hP_{X^{t,\xi}_r})dr
+ \int_t^s \si(r, X^{t,\xi}_r, \hP_{X^{t,\xi}_r}) dB_r
+\b_s + K^{t,\xi}_s, \qq s\in[t, T], \nonumber 
\eea
and then using the law $\hP_{X^{t,\xi}}$ to consider a slight variation of (\ref{Xtxi}):
\bea
\label{Xtxxi}
X_{s}^{t,x,\xi} &=& x + \int_t^s \int_{A \times \hR_{+}}  \th(r, X^{t,x,\xi}_{r-}, \hP_{X^{t,\xi}_{r}},  z) 
\1_{ [0, \l(r, X^{t,x,\xi}_{r-}, \hP_{X^{t,\xi}_{r}})] } (y) \wt\cN^{s}(dr dz dy)\\
&& +\int_t^s b(r, X^{t,x,\xi}_r, \hP_{X^{t,\xi}_r})dr +\int_t^s \si(r, X^{t,x,\xi}_r, \hP_{X^{t,\xi}_r})dB_r +\b_s + K^{t,x,\xi}_s, \qq s\in[t, T], \nonumber
\eea
where  $\xi \in L^{2}(\mathcal{F}_{t}; \mathbb{R})$,
then we shall argue below that the following flow property holds:
\bea
\label{flow}
\left(X_{r}^{s,X_{s}^{t,x,\xi}, X_{s}^{t,\xi}}, X_{r}^{s, X_{s}^{t,\xi}} \right) 
= (X_{r}^{t,x,\xi}, X_{r}^{t,\xi} ), \qq 0 \leq t \leq s \le r \leq T, 
\eea
for all $ (x,\xi) \in \mathbb{R}\times L^{2}(\mathcal{F}_{t}; \mathbb{R})$. We should note that although both SDEs (\ref{Xtxi}) and (\ref{Xtxxi}) resemble the
original equation (\ref{MFSDER}), the process $X^{t,x,\xi}$ has the full information of the solution given the initial data $(x, \xi)$, where
$\xi$ provides the initial distribution $\hP_{\xi}$, and $x$ is the actual initial state. 

To prove the well-posedness of SDEs (\ref{Xtxi}) and (\ref{Xtxxi}), we first recall the so-called ``Discontinuous Skorohod 
Problem" (DSP) (see, e.g., \cite{DI91, Ma93}). Let $Y\in \hD([0,T])$, $Y_0\ge0$. We say that a pair $(X, K)\in \hD([0,T])^2$ is a solution to
the DSP($Y$) if 

\ss
(i) $X=Y+K$; 

\ss

(ii) $X_t\ge 0$, $t\ge 0$; and 

\ss

(iii) $K$ is nondecreasing, $K_0=0$, and $K_t=\int_0^t \1_{\{X_{s-}=0\}}dK_s$, $t\ge 0$.

It is well-known that the solution to DSP exists and is unique, and it can be shown (see \cite{Ma93}) that the condition (iii) amounts to saying that $\int_0^t\1_{\{X_{s-}>0\}}dK^c_s=0$, where $K^c$ denotes the continuous part of $K$, and  $\D K_t= (X_{t-}+\D Y_t)^-$. Furthermore, it is shown in \cite{DI91} that solution mapping of the 
DSP, $\Gamma: \hD([0,T]) \mapsto \hD([0,T])$, defined by $\Gamma(Y) =X$, is Lipschitz continuous under uniform topology. That is,
 there exists a constant $L>0$ such that 
\bea
\label{GammaLip}
\sup_{t \in [0,T]} |\Gamma(Y^1)_t- \Gamma(Y^2)_t|  \leq  L \sup_{t\in [0,T]} |Y_{t}^{1} - Y_{t}^{2}|, \qq Y^1, Y^2\in \hD([0,T]). 
\eea

Before we proceed to prove the well-posedness of (\ref{Xtxi}) and (\ref{Xtxxi}), we note that the two SDEs can be argued separately. 
Moreover, while (\ref{Xtxi}) is a mean-field (or McKean-Vlasov)-type of SDE, (\ref{Xtxxi}) is actually a standard SDE (although with 
state-dependent intensity) with discontinuous paths and reflection, given the law of the solution to (\ref{Xtxi}), $\hP_{X^{t,\xi}}$, 
and it can be argued similarly but much simpler. 
Therefore, in what follows we shall focus only on the well-posedness of SDE (\ref{Xtxi}). Furthermore, for simplicity we shall assume $b\equiv 0$, as the general case can be
argued similarly without substantial difficulty. 

The scheme of solving the SDE (\ref{Xtxi}) is more or less standard (see, e.g., \cite{Ma93}). We shall first consider an SDE without
reflection: for  $\xi\in L^2(\cF_t; \hR)$ and $s\in[t, T]$, 
\bea
\label{Ytxi}
Y_{s}^{t,\xi} &=&
\xi + \int_t^s \int_{A \times \dbR_{+}}  \th(r, \Gamma(Y^{t,\xi})_{r-},  \hP_{\G(Y^{t,\xi})_{r}},  z) \1_{ [0,\l_{r-}^{\Gamma(t,\xi)} ] } (y) \wt\cN^{s}(dr dz dy) \\
&& + \int_t^s \si (r, \Gamma(Y^{t,\xi})_{r},  \hP_{\G(Y^{t,\xi})_{r}}) dB_r + \beta_s, \nonumber
\eea
where $\lambda_{r-}^{\Gamma(t,\xi)} := \l(r, \G(Y^{t,\xi})_{r-}, \hP_{\G(Y^{t,\xi})_{r}})$. Clearly, if we can show that (\ref{Ytxi}) 
is well-posed, then by simply setting  $X^{t,\xi}_s = \G(Y^{t,\xi})_s$ and $K_{s}^{t, \xi}= X^{t, \xi}_s-Y^{t, \xi}_s$, $s\in[t,T]$, 
we see that $(X^{t,\xi}, K^{t, \xi})$ would solve  SDE (\ref{Xtxi})(!). We should note that a technical difficulty caused by the
 presence of the state-dependent intensity is that the usual $L^2$-norm does not work as naturally as expected, as we shall
 see below. We nevertheless have the following result.  


\begin{thm}
\label{thm1-MFSDER}
Assume that Assumptions \ref{assump1} is in force. Then, 
there exists a solution $Y^{t,\xi}\in L^2_\hF(\hD([t,T]))$ to SDE (\ref{Ytxi}).
Furthermore, such solution is pathwisely unique. 
\end{thm}

{\it Proof.} 
Assume $t=0$. For a given $T_0>0$, and $y\in L^1_\hF (\hD([0, T_0]) )$, consider a mapping $\sT$:
\bea
\sT(y)_s &:=& \xi 
+ \int_0^{s} \int_{A \times \dbR_{+}} \th(r, \G(y)_{r-}, \hP_{\G(y)_r}, z) \1_{[0, \lambda(r, \Gamma(y)_{r-}, \hP_{\Gamma(y)_r})] } (u) \wt\cN^s (drdz du) \\
&& + \int_0^s \si (r, \G(y)_r, \hP_{\G(y)_r}) dB_r + \beta_s, \qq s\ge 0. \nonumber
\eea
We shall argue that $\sT$ is a contraction mapping on $L^1_\hF (\hD([0, T_0]))$ for $T_0 >0$ small enough.  

To see this, denote, for $\eta\in\hD([0,T_0])$, $|\eta|^{*}_s:=\sup_{0\le r\le s}|\eta_r|$, and   
define $\th_s(z) := \th(s, \Gamma(y)_{s},\hP_{\Gamma(y)_{s}}, z)  $, $\l_s:= {\lambda} (s, \Gamma(y)_{s}, \hP_{\Gamma(y)_{s}})$, $\si_s:= \si(s, \Gamma(y)_{s},\hP_{\Gamma(y)_{s}} )$, $s\in[0,T_0]$. Then, we have
\beaa
\dbE[ |\sT(y) |_{T_0}^*] \neg&\neg \leq\neg & \neg C \Big\{\hE|\xi|+\hE \Big[ 
\int_0^{T_0}\neg \int_{A \times \dbR_{+}}
\big|\th_r(z)\1_{ [0,\lambda_r]} (y))\big|  \n^s (dz)dy dr \Big]  + \hE\Big[\Big(\int_0^{T_0}  | \si_r |^2 dr\Big)^{1/2}\Big] \Big\}  \\
\neg&\neg \leq\neg & \neg C \hE|\xi|
+ C \hE \Big[ \int_0^{T_0} \int_{A} | \th_r(z) \lambda_r  | \nu^s(dz) dr \Big]   + C \hE\Big[\Big(\int_0^{T_0}  | \si_r |^2 dr\Big)^{1/2} \Big] < \infty,
\eeaa
thanks to Assumption \ref{assump1}. Hence, $\sT(y)\in  L^1_\hF(\hD([0, T_0]) ) $.

We now show that  $\sT$ is a contraction on $L^1_\hF(\hD([0, T_0]) )$. For $y_1, y_2 \in L^1_\hF(\hD([0, T_0]))$, we denote  
$\th^i$, $\l^i$, and $\si^i$, respectively, as before, and denote $\Delta \f  := \f ^1 -\f ^2$, for $\f=\th, \l, \si$, and $\D\sT(s) =\sT(y_1)_s-\sT(y_2)_s$, 
$s\ge 0$.  
Then, we have, for $s \in [0,T_0]$,
\beaa
\D\sT(s)= 
\int_0^{s} \int_{A \times \dbR_{+}} \left[\D\th_r(z) \1_{ [0,\lambda_r^1] } (y)+\th_r^2 (z)(\1_{ [0,\lambda_r^1] } (y)-\1_{ [0,\lambda_r^2] } (y))\right]\wt\cN^{s}(dr dz dy) + \int_0^s \Delta \si_r dB_r.
\eeaa
Clearly, $\D\sT=\sT(y_1) - \sT(y_2)$ is a martingale on $[0, T_0]$.  Since $\wt \cN=\cN-\h \cN$ and $|\1_{[0,a]}(\cd)-\1_{[0,b]}(\cd)|\le \1_{[a\wedge b, a\vee b]}(\cd)$ for any $a,b\in\hR$, we have, for $0\le s\le T_0$, 
\bea
\label{T12}
\hE | \D\sT |^*_s  \neg&\neg\le\neg&\neg 2 \hE \Big[ 
\int_0^s\neg\neg\int_{A \times \dbR_{+}} \neg
\big| \D\th_r(z)\1_{ [0,\lambda_r^1] } (y)+\th_r^2(z) (\1_{ [0,\lambda_r^1] } (y)-\1_{ [0,\lambda_r^2] } (y))\big|  \n^s (dz)dy dr \Big] 
\nonumber \\
&&+ \hE \Big[ \Big( \int_0^s |\D \si_r|^2 dr \Big)^{\frac{1}{2}} \Big] 
:= I_1+I_2. 
\eea
 Recalling from Remark \ref{Remark1_S2}-(ii) for the generic constant $C>0$, and by Assumption \ref{assump1}-(iv), (\ref{GammaLip}),  and the definition of $W_1(\cd, \cd)$ (see (\ref{W2})), we have
\bea
\label{I12}
I_1 &\le & C \Big\{\hE\Big[\int_0^s\neg \int_{A}| \D\th_r(z)| \n^s (dz)dr\Big]+\hE\Big[\int_0^s |\D\l_r | dr \Big]\Big\} \nonumber \\
&\le & C \hE\Big[\int_0^s [|\G(y_1)_r-\G(y_2)_r|  + W_1 (\hP_{\Gamma(y_1)_r}, \hP_{\Gamma(y_2)_r}) ] dr \Big] \nonumber \\
&\le & C \hE\Big[\int_0^s |\G(y_1)-\G(y_2)|^{*}_r dr\Big] \le C T_0 \| y_1 - y_2 \|_{L^1(\hD([0,T_0]))}, \\
I_2 &\le & C \hE \Big[ \Big( \int_0^s \{ | y_1-y_2|^{*, 2}_r + W_1 (\hP_{\G(y_1)_r}, \hP_{\G(y_2)_r})^2 \} dr \Big)^{1/2}\Big] \nonumber\\ 
&\le &  C \hE \big[ \sqrt{s} \big( |y_1-y_2|_s^{*} + \hE |y_1-y_2|_s^{*} \big) \big]   
\le C \sqrt{T_0} \| y_1-y_2\|_{L^1(\hD([0,T_0]))}. \nonumber
\eea
Combining (\ref{T12}) and (\ref{I12}), we deduce that 
\bea
\label{DYI1}
\| \D\sT\|_{L^1(\hD([0,T_0]))} \leq C(T_0 + \sqrt{T_0})  \| y_1-y_2\|_{L^1(\hD([0,T_0]))}, \qq s \in [0, T_0].
\eea
Therefore, by choosing $T_0$ such that $C(T_0 + \sqrt{T_0}) < 1$, we see that the mapping $\sT$ is a contraction on
$L^1(\hD([0,T_0])) $, which implies that (\ref{Ytxi}) has a unique solution in $L^1_\hF(\hD([0, T_0]) )$.  Moreover, 
we note that $T_0$ depends only 
on the universal constants in Assumption \ref{assump1}. We can repeat the argument for the time interval $[T_0, 2 T_0], [2T_0, 3T_0], \cdots$, and conclude that  (\ref{Ytxi}) has a unique solution in $L^1_\hF (\hD([0, T]) )$ for any given $T>0$. 

Finally, we claim that the solution $Y \in L^2_\hF (\hD([0,T]))$. Indeed, by Burkholder-Davis-Gundy's inequality and Assumption \ref{assump1},
we have
\bea
\label{Yestp}
\hE [  |Y|^{*,2}_s]
\neg&\neg \leq\neg & \neg C \Big\{\hE|\xi|^2+\hE \Big[ 
\int_0^s\neg \int_{A \times \dbR_{+}}
\big|\th_r( z)\1_{ [0,\lambda_r]} (y))\big|^2  \n^s (dz)dy dr\Big] + \hE[ \int_0^s  | \si_r|   ^2 dr]
+\hE|\b|^{*,2}_T\Big\}\nonumber \\
\neg&\neg \leq\neg & \neg C\Big\{\hE|\xi|^2+ \hE \Big[ 
\int_0^s\neg\big[1+|Y_r|^2+W_1(0, \G(Y)_r)\big]^2  dr \Big]+\hE|\b|^{*,2}_T\Big\}\\
\neg&\neg \leq\neg & \neg C\Big\{\hE|\xi|^2+\int_0^s(1+ \hE[|Y|^{*,2}_r])dr+\hE|\b|^{*,2}_T\Big\}, \qq s\in[0,T]. \nonumber
\eea
Here, in the last inequality above we used the fact that 
$$W_1(0, \G(Y)_r)^2\le (\|\G(Y)_r\|_{L^1(\O)})^2\le (\hE|\G(Y)|^*_r)^2\le C \hE[|Y|^{*,2}_r], \qq r\in[0,s].
$$
 Applying the Gronwall inequality, we obtain that $\hE [  |Y|^{*,2}_T] <\infty$. The proof is now complete. 
\qed


\begin{rem}
{\rm
(i) It is worth noting that once we solved $X^{t,\xi}$, then we know $\hP_{X^{t,\xi}}$, and (\ref{Xtxxi}) can be viewed as a  standard SDEDR with coefficient $\tilde{\lambda}(s, x) := {\lambda}(s,x, \hP_{X_{s}^{t,\xi}}) $, which is Lipschitz in $x$. This guarantees the existence and uniqueness
of the solution $(X^{t,x,\xi}, K^{t,x,\xi})$ to (\ref{Xtxxi}).

\ms
(ii) The uniqueness 
of the solutions to (\ref{Xtxi}) and (\ref{Xtxxi}) implies that  $X_{s}^{t,x, \xi}\mid_{x=\xi} = X_{s}^{t,\xi} , \; s \in [t,T]$.  That is, $X_{s}^{t,x, \xi}\mid_{x=\xi}$ solves the same SDE as $X_s^{t,\xi}, \; s \in [t,T]$. (See more detail in \cite{thesis}) 

\ms
(iii) Given $(t,x) \in [0,T] \times \mathbb{R}$, if $\hP_{\xi_1} = \hP_{\xi_2}$ for $\xi_1, \xi_2 \in L^2(\mathcal{F}_t; \mathbb{R})$, then $X^{t,x,\xi_1}$ and $X^{t,x,\xi_2}$ are indistinguishable. So,
$X^{t,x,\hP_{\xi}} := X^{t,x,\xi}$, i.e. $X^{t,x,\xi}$ depends on $\xi$ only through its law.  
\qed 
}
\end{rem}

\subsection{An It\^o's formula}


%
%
%

We shall now present an It\^o's formula that will be frequently used in our future discussion. We note that a similar formula
for mean-field SDE can be found in \cite{BLPR}, and the one involving jumps was given in the recent work \cite{MFjump}. The one 
presented below is a slight modification of that of \cite{MFjump}, taking the particular  state-dependent intensity feature of the dynamics
into account. Since the proof is more or less standard  but quite lengthy, we refer to \cite{thesis} for the details. 

In what follows we let $(\tilde{\Omega}, \tilde{\mathcal{F}}, \tilde{\hP})$ be a copy of the probability space $(\Omega, \mathcal{F}, \hP)$,
and denote $\tilde{\dbE}[\cdot]$ to be the expectation under $\tilde{\hP}$. For any random variable $\vartheta$ defined on $(\Omega, \mathcal{F}, \hP)$, we denote, when there is no danger of confusion,  $\tilde{\vartheta} \in (\tilde{\Omega}, \tilde{\mathcal{F}}, \tilde{\hP})$ to be a copy of $\vartheta$ such that $\tilde{\hP}_{\tilde{\vartheta}} = \hP_{\vartheta}$. We note that that  $\tilde{\dbE}[\cdot]$ acts only on the variables of the form $\tilde\vartheta$. 

%
%

We first define the following classes of functions.
\begin{defn}
\label{def_C^1,2,1}
We say that $F \in C_{b}^{1,2,(2,1)} ([0,T] \times \mathbb{R}  \times \mathbb{R} \times \sP_{2}(\mathbb{R}) )$, if 

(i) $F(t,v,\cdot, \cdot) \in C_{b}^{2,1} (\mathbb{R} \times \sP_{2}(\mathbb{R}) )$, for all $t \in [0,T]$ and $v \in \hR$;
\ss

(ii) $F(\cdot, v,x,\mu) \in C_{b}^{1}([0,T])$, for all $(v, x,\mu) \in  \hR \times \mathbb{R} \times \sP_{2}(\mathbb{R})$; 
\ss

(iii) $F(t, \cdot, x, \mu) \in C_b^{2} (\mathbb{R})$, for all $(t,x,\mu) \in [0,T] \times \mathbb{R} \times \sP_{2}(\mathbb{R}) $;

(iv) All derivatives involved in the definitions above are uniformly bounded over $[0,T] \times \mathbb{R} \times \mathbb{R} \times \sP_{2}(\mathbb{R}) $ and Lipschitz in $(x,\mu)$, uniformly with respect to $t$. 
\qed
\end{defn}

We are now ready to state the It\^o's formula. 
Let $V^{t,v}$ be an It\^o process given by
\bea
\label{SDE_V}
V_s^{t,v} = v + \int_t^s b^{V}(r,V_r^{t,v}) dr + \int_t^s \sigma^V (r, V_r^{t,v}) dB^V_r
\eea
where $v\in \hR$ and $(B^V_t)_{t \in [0,T]}$ is a standard Brownian motion independent of $(B_t)_{t \in [0,T]}$.
For notational simplicity, in what follows for the coefficients $\varphi=b,\si, \b, \l$, we 
denote $\f^{t, x,\xi}_s:=\f(s, X^{t,x,\xi}_s, \hP_{X^{t, \xi}_s})$,
 $ \th^{t,x,\xi}_s(z):=\th(s, X^{t,x,\xi}_s, \hP_{X^{t, \xi}_s},z)$,
 $\tilde{\f}^{t,\xi}_s:=\f(s, \tilde{X}^{t,\tilde{\xi}}_s, \hP_{X^{t, \xi}_s})$,
 and
 $ \tilde{\th}^{t,\xi}_s(z):=\th(s, \tilde{X}^{t,\tilde{\xi}}_s,$ $ \hP_{X^{t, \xi}_s},z)$.
 Similarly, denote $b_s^{t,v} := b^V (s, V_s^{t,v})$ and $\si_s^{t,v} := \si^V (s, V_s^{t,v})$.
 Also, let us write $\Theta_s^t := (s, V_s^{t,v}, X_{s}^{t,x,\xi}, \hP_{X_{s}^{t,\xi}})$. Then $\Th^t_t=(t,v,x,\hP_\xi)$.
\begin{prop}[It\^o's Formula]
\label{Ito}
Let $\Phi \in C_{b}^{1,2,(2,1)} ([0,T] \times \mathbb{R}\times \mathbb{R} \times \sP_{2}(\mathbb{R}) )$, and
 $(X^{t, \xi}, X^{t,x,\xi}, V^{t,v})$ be the solutions to 
(\ref{Xtxi}), (\ref{Xtxxi}) and (\ref{SDE_V}), respectively,  
on $[t,T]$.  
Then, for $0 \leq t \leq s \leq T$, it holds that
\bea
\label{ItoForm}
&& \Phi(\Th^t_s) - \Phi(\Th^t_t) \nonumber\\
\neg&\neg = \neg&\neg \int_{t}^{s}  \Big( 
\partial_{t} \Phi(\Theta_r^t) + \partial_{x} \Phi(\Theta_r^t) b^{t,x,\xi}_r + \frac{1}{2} \partial_{xx}^{2}\Phi(\Theta_r^t) (\si^{t,x,\xi}_r)^2 
+ \partial_{v} \Phi(\Theta_r^t) b_r^{t,v} + \frac{1}{2} \partial_{vv}^{2}\Phi(\Theta_r^t) (\si_r^{t,v})^2  
 \Big) dr \nonumber\\
&&\neg\neg+ \int_{t}^{s} \partial_{x}\Phi (\Theta_r^t) \si^{t,x,\xi}_{r} dB_r 
+  \int_{t}^{s} \partial_{v}\Phi (\Theta_r^t) \si_r^{t,v} dB^V_r 
+ \int_t^s \partial_{x}\Phi (\Theta_{r-}^t) \1_{ \{ X_{r-}=0 \} } dK_r \\
&&\neg\neg  +
\int_{t}^{s} \int_{A} \Big( \Phi(r, V_r^{t,v}, X_{r-}^{t,x,\xi} + \th^{t,x,\xi}_{r}( z), \hP_{X_{r}^{t,\xi}}) - \Phi (\Th_{r-}^t)  - \partial_{x} \Phi(\Th^t_{r-}) \th^{t,x,\xi}_{r-}( z) \Big) \lambda^{t,x,\xi}_{r} \nu^s(dz) dr \nonumber\\
&&\neg\neg+\int_{t}^{s} \int_{A \times \dbR_{+}} \left\lbrace 
\Phi(r, V_r^{t,v}, X_{r-}^{t,x,\xi} + \th^{t,x, \xi}_{r}(z), \hP_{X_{r}^{t,\xi}}) 
- \Phi(\Th^t_{r-}) \right\rbrace 
\1_{[0, \lambda^{t,x,\xi}_r]} (y) \wt\cN^s(dr dz dy) \nonumber\\
%
%
%
&&\neg\neg+  \int_{t}^{s} \wt{\mathbb{E}} \Big[ \partial_{\mu}\Phi (\Theta_r^t,  \tilde{X}_{r}^{t,\tilde{\xi}}) \tilde b^{t, {\xi}}_{r}
+ \frac{1}{2} \partial_{y}(\partial_{\mu}\Phi) (\Theta_r^t,  \tilde{X}_{r}^{t,\tilde{\xi}}) (\tilde \si^{t, {\xi}}_{r})^2 \nonumber \\
%
%
&&\qq \q + \int_{0}^{1}\neg\neg \int_{A} [\partial_{\mu} \Phi (\Theta_r^t, \tilde{X}_{r}^{t,\tilde{\xi}}\neg +\neg \rho \tilde\th^{t,x, \xi}_{r}( z)) 
\neg-\neg \partial_{\mu}\Phi (\Theta_r^t, \tilde{X}_{s}^{t,\tilde{\xi}}) ] 
\tilde \th^{t,{\xi}}_r( z) \tilde\lambda^{t,{\xi}}_{r}\nu^s(dz) d\rho \Big]dr. \nonumber
\eea
\qed
\end{prop}

\section{A Bertrand game among the sellers (static case)}

In this section we analyze a price setting mechanism among liquidity providers (investors placing sell limit orders), which will be used as the 
basis for our continuous time model in the rest of the paper.  Following the ideas of \cite{CS15, LS11, LS12}, we shall consider this process of the (static) price setting as a Bertrand-type of game among the sellers, 
each placing a certain number of sell limit orders at a specific price, and trying to maximize her expected utility. 
To be more precise, we assume that 
sellers use the price at which they place limit orders as their {\it strategic variable}, 
and the number of shares submitted would be determined accordingly. Furthermore, we assume that 
there is a {\it waiting cost}, also as a function of the price. Intuitively, a higher price will lead to a longer execution 
time, hence a higher waiting cost. Thus, there is a competitive game among the sellers for better total reward. 
Finally, 
we assume that the sellers are homogeneous in the sense that they have the same subjective probability measure, 
so that they share the same degree of risk aversion (or uncertainty aversion). 


We now give a brief description of the problem. We assume that there are $N$ sellers, and the $j$th seller places limit orders at price 
$p_j = X + l_j$, $j=1,2,\cdots,N$, where $X$ is the mid price. Without loss of generality, we may assume $X=0$. 
As a main element in an {\it oligopolistic competitions} (cf. e.g., \cite{LS12}), we assume that there exists a {\it demand function}, 
denoted by $h_i^N (p_1, p_2, \cdots, p_N)$ for each seller $i$, at a  given price vector $p=(p_1, p_2, \cdots, p_N)$ at the moment. We note that the number of shares of limit orders the seller $i$ places in the LOB will be determined by the value of her demand function with the given price vector, 
{hence a {\it Bertrand game}.}\footnote{A {\it Cournot game} is one such that the price $p_i$ is the function of 
the numbers of shares $q=(q_1, \cds, q_N)$ through a demand function. The two games are often exchangeable if the 
demand functions are invertible (see, e.g., \cite{LS12}).} More specifically, we assume that the demand functions $h_i^N$, $i=1,2,\cdots, N$,  are smooth  and satisfy the following properties:
\bea
\label{mono}
\frac{\partial h_i^N}{\partial p_i} < 0, \hspace{0.3cm} \text{ and }
\hspace{0.3cm}  \frac{\partial h_i^N}{\partial p_j} > 0, \text{~~ for } j \neq i.
\eea
We note that (\ref{mono}) amounts to saying that the number of shares each seller places is decreasing in the seller's own price and increasing in the other sellers' price. Furthermore, we shall assume that the demand functions are invariant under permutations of the other sellers' prices, in the sense that, for fixed $p_1, \cdots, p_N$ and all $i,j \in \{ 1,\cdots, N \}$,
\bea
\label{exchange}
h_i^N (p_1, \cdots, p_i, \cdots, p_j, \cdots, p_N)= h_j^N (p_1, \cdots, p_j, \cdots, p_i, \cdots, p_N). 
\eea
It is worth noting that the combination of (\ref{mono}) and (\ref{exchange}) is the following fact: if a price vector $p$ is ordered by
$p_1 \leq p_2 \leq \cdots \leq p_N$, then by (\ref{mono}) and (\ref{exchange}), for any $i < j$, it holds that
\bea
\label{order-h}
h_j^N (p)  = h_i^N ( p_1, \cdots, p_j, \cdots, p_i, \cdots, p_N) \le h_i^N (p_1, \cdots, p_i, \cdots, p_i, \cdots, p_N) \le h_i^N (p). 
\eea
That is, the demand functions are ordered in a reversed way. 
%
Finally, for each $i$, we denote the price vector for ``other" prices for seller $i$ by $p_{-i}$,  and we 
assume that there is a ``least favorable" price for seller $i$, denoted by $\hat{p}_i (p_{-i}) < \infty$, in the sense that
\bea
\label{Choke}
h_i^N (p_1,\cdots,p_{i-1},\hat{p}_i, p_{i+1},\cdots,p_N) = 0.
\eea
The price $\hat p_i$ is often called  the 
{\it choke price}. We note that the existence of such price, together with the monotonicity property (\ref{mono}), indicates the possibility that $h^N_j(p)<0$, for
some $j$ and some price vector $p$. But, on the other hand, since the size of order placement cannot be negative, such scenario 
becomes unpractical. To amend this, we introduce the notion of \textit{actual demand}, denoted by $\{\h h_i(p)\}$,  which we now describe. 

Consider an ordered price vector $p=(p_1, \cds, p_N)$, with $p_i\le p_j$, $i\le j$, and we look at $h_N^N (p)$.  
If $h_N^N (p) \geq 0$, then by (\ref{order-h})  we have $h_i^N (p) \geq 0$ for all $i = 1,\cdots, N$. In this case, we denote $\h h_i(p)= h_i^N (p)$ for all $i=1, \cdots, N$.
If $h_N^N (p) < 0$, then we set $\h h_N (p) = 0$. That is, the $N$-th seller does not act at all. We assume that the 
remaining $N-1$ sellers will observe this fact and modify their strategy as if there are only $N-1$ sellers. More precisely, we first choose a choke price $\h p_N$ so that
$h^N_N(p_1, \cds, p_{N-1}, \h p_N)=0$, and define
$$ h_i^{N-1} (p_1, p_2, \cdots, p_{N-1}) := h_i^N (p_1, p_2, \cdots, p_{N-1}, \h{p}_N), \qq  i=1,\cdots, N-1,
$$ 
and continue the game among the $N-1$ sellers.

 In general, for $1\le n\le N-1$, assume the $(n+1)$-th demand functions $\{h^{n+1}_i\}_{i=1}^{n+1}$ are defined. If $h_{n+1}^{n+1}(p_1, \cds, p_{n+1}) < 0$, then other $n$ sellers will assume ($n+1$)-th seller sets a price at $\hat{p}_{n+1}$ with zero demand (i.e., $h_{n+1}^{n+1} (p_1, p_2, \cdots, p_n, \hat{p}_{n+1}) = 0$), and modify their demand functions to 
\bea
\label{hn}
 h_i^n (p_1, p_2, \cdots, p_n) :=  h_i^{n+1} (p_1, p_2, \cdots, p_n, \hat{p}_{n+1}), \qq i=1,\cdots, n.
 \eea
We can now define the ``actual demand function" $\{\h h_i\}_{i=1}^N$. 
%
%
\begin{definition}[{{\it Actual demand function}}]
\label{ActualDemand}
Assume that $\{h^N_i\}_{i=1}^N$ is a family of demand functions. The family of  ``actual demand functions", denoted by 
$\{\h h_i\}_{i=1}^N$, 
 are defined in the following
steps: for a given ordered price vector $p$,

(i)  if $h_N^N(p) \geq 0$, then 
we set $\h h_i (p) = h_i^N (p)$ for all $i= 1,\cdots,N$;

(ii) if $h^N_N(p)<0$, then we define recursively for  $n=N-1, \cds 1$ the demand functions $\{h^n_i\}_{i=1}^n$ as in (\ref{hn}). In particular, if there exists an $n< N$ such that $ h_{n+1}^{n+1} (p_1, p_2, \cdots, p_n, {p}_{n+1}) < 0$ and
$h_n^n (p_1, p_2, \cdots, p_n) \geq 0$, then we set 
\bea
\label{hathn}
\h h_i (p) = 
\begin{cases}
h_i^n ( p_1, p_2, \cdots, p_n) & i =1, \cdots,n \\ 
0 &  i = n+1, \cdots, N; 
\end{cases} 
\eea

(iii) if there is no such $n$, then $\h h_i (p) = 0$ for all $i =1, \cdots, N$.
\qed
\end{definition}

We note that the actual demand function will always be non-negative, but for each price vector $p$, the number $\#\{i : \h h_i(p)>0\}\le N$,
and could even be zero.

\subsection{The Bertrand game and its Nash equilibrium}

Besides the demand function, a key ingredient in the placement decision making process is the 	``waiting cost" for the time it takes
for the limit order to be executed. We shall assume that each seller has her own waiting cost function $c_i^N \triangleq c_i^N (p_1, p_2, \cdots, p_N,Q)$, where $Q$ is the total number of shares available in the LOB. Similar to the demand function, we shall assume the following assumptions for the waiting cost. 
\begin{assum}
\label{Assump_cost1}
For each seller $i \in \{1, \cdots, n \}$ with $n \in [1, N]$,  each $c_i^N$ is smooth in all variables  such that

(i) (Monotonicity)  $\frac{\partial c_i^N}{\partial p_i} > 0$, and  $\frac{\partial c_i^N}{\partial p_j} < 0$, for $j \neq i$;

\ms
(ii) (Exchangeability) 
$c_i^N (p_1, \cdots, p_i, \cdots, p_j, \cdots, p_N)
= c_j^N (p_1, \cdots, p_j, \cdots, p_i, \cdots, p_N)$; 

\ms
(iii) $c^N_i(p)\big|_{p_i=0}=0$, and $\frac{\partial c_i^N}{\partial p_i} \Big|_{p_i=0+}\in (0, 1)$;

\ms
(iv) $\lim_{p_i\to \infty}\frac{p_i}{c^N_i(p)}=0$, $i=1,\cds, N$.
\qed
\end{assum}

\begin{rem}
{\rm
(i) By exchangeability, in what follows, we shall assume without loss of generality that all prices are ordered.

(ii) Assumption \ref{Assump_cost1}-(i), (ii) ensure that the price ordering leads to the same ordering for waiting cost functions, similar to 
what  we argued before for demand functions.
%

(iii) Consider the function $J_i(p, Q)=p_i-c^N_i(p, Q)$. Assumption \ref{Assump_cost1} amounts to saying that $J_i(p, Q)\big|_{p_i=0}=0$, 
$\frac{\pa J_i(p,Q)}{\pa p_i}\Big|_{p_i=0+}>0$, and $\lim_{p_i\to \infty}J_i(p, Q)<0$. Thus, there exists $p_i^0=p^0_i(p_{-i}, Q)>0$ such that
$\frac{\pa J_i(p, Q)}{\pa p_i}\Big|_{p_i=p^0_i} =0$, and  $\frac{\pa J_i(p, Q)}{\pa p_i}\Big|_{p_i>p^0_i} <0$. 

(iv) Since $J_i(0,Q)=0$, and $\frac{\pa J_i(p,Q)}{\pa p_i}\Big|_{p_i=0+}>0$, one can easily check that $J_i(p^0_i, Q)>0$. This, 
together with Assumption \ref{Assump_cost1}-(iv), shows that there exists $\tilde{p}_i=\tilde p_i(p_{-i},Q)>p^0_i$, such that $J_i(p_i, Q)\big|_{p_i=\tilde p_i}=0$ (or, equivalently $c_i^N(p_1, \cdots, p_{i-1}, \tilde{p}_i, p_{i+1}, \cdots, p_N, Q)=\tilde{p}_i$). Furthermore, the Remark (iii) 
implies that $J_i(p_i, Q)<0$ for all $p_i>\tilde p_i(p_{-i},Q)$.
In other words, any selling price higher than $\tilde p_i(p_{-i},Q)$ would yield a negative profit, and therefore should be prevented. 
\qed}
\end{rem}

The Bertrand game among sellers can now be formally introduced: each seller chooses its price to maximize profit in a non-cooperative manner, and their decision will be based not only on her own price, but also on the actions of all other sellers.
We denote the profit of each seller by 
\bea
\label{profit}
\Pi_{i} (p_1, p_2, \cdots, p_N, Q) := \h h_{i}(p_1, p_2, \cdots, p_N)[p_i - c_i^N(p_1, p_2, \cdots, p_N, Q)],
\eea
and each seller tries to maximize her profit $\Pi$. 
For each fixed $Q$, we are looking for  
a Nash equilibrium price vector 
$p^{*,N}(Q)=(p_1^{*,N}(Q), \cdots, p_N^{*,N}(Q))$. We note that in the case when 
$\h h_i(p^{*,N})=0$ for some $i$, the $i$-th seller will not participate in the game (with zero profit), so we shall modify the price 
\begin{equation}
p^{*,N}_i(Q) \triangleq c_i^N (p_1^{*,N}, \cdots, p_N^{*,N}, Q) = c_i^N (p^{*,N}, Q),
\end{equation}
and consider a subgame involving the $N-1$ sellers, and so on. 
That is, for a subgame with $n$ sellers, they solve
\begin{equation}
\label{subgame_NE}
p_i^{*,n} = \arg \max_{p \geq 0} \Pi_i^n (p_1^{*,n}, p_2^{*,n}, \cdots, p_{i-1}^{*,n}, p, p_{i+1}^{*,n}, \cdots, p_n^{*,n}, Q),
\hspace{0.2cm} i=1,\cdots, n
\end{equation}
to get $p^{*,n} = (p_1^{*,n}, \cdots, p_n^{*,n}, c_{n+1}^{*,n+1}, \cdots, c_{N}^{*,N})$.
More precisely, we define a Nash Equilibrium as follows. 
\begin{definition}
A vector of prices $p^*=p^{*} (Q)= (p_1^*, p_2^*, \cdots, p_N^*)$ is called a Nash equilibrium if
\bea
\label{equil}
p_i^*= \arg \max_{p \geq c_i } \Pi_i (p_1^*, p_2^*, \cdots, p_{i-1}^*, p, p_{i+1}^*, \cdots, p_N^*, Q),
\eea
and  $p_i^* = c_i^{*,i} (p^*,Q) $ whenever $\h h_i(p^*) = 0$, $i=1,2, \cdots, N$.
\qed
\end{definition}

We assume the following on a subgame for our discussion.

\begin{assum}
For $n=1, \cdots, N$, we assume that there exists a unique solution to the system of maximization problems in equation (\ref{subgame_NE}).
\end{assum}
%
%
%
\begin{rem}
\label{Nash}
{\rm We observe from Definition of the Nash Equilibrium that, in equilibrium, a seller is actually participating in the Bertrand game only 
when
her actual demand function is positive, and those with zero actual demand function will be ignored in the subsequent subgames. However, a participating seller does not necessarily have positive profit unless she sets the price higher than the waiting cost. In other
words, it is possible that $\h h_i(p^*)>0$, but $p^*_i=c_i(p^*,Q)$, so that $\Pi_i(p^*,Q)=0$. We refer to such a case the {\it boundary}
case, and denote the price to be $c_i^{*,b}$. 
\qed
}
\end{rem}

The following result details the procedure of finding the Nash equilibrium for the Bertrand competition. The idea is quite similar to 
that in \cite{CS15}, except for the general form of the waiting cost. We sketch the proof for completeness.
\begin{prop}
\label{existence1}
Assume that Assumption \ref{Assump_cost1} is in force. Then there exists a 
Nash equilibrium to the Bertrand game 
(\ref{profit}) and (\ref{equil}). 

Moreover, the equilibrium point $p^*$, after the modifications, should take the following form:
\bea
\label{equilip}
p^* = ( p_1^{*} , \cdots, p_k^{*}, c_{k+1}^{*, b}, \cdots, c_n^{*,b}, c_{n+1}^{*}, \cdots, c_N^*),
\eea
from which we can immediately read: $\h h_i(p^*)>0$ and $p_i^*> c^*_i$,  $i=1, \cds, k$; $\h h_i(p^*)>0$ but $p_i^*\le c^*_i$,  $i=k+1, \cds, n$; and $\h h_i(p^*)\le 0$, $i=n+1, \cds, N$.\end{prop}

{\it Proof.}
We start with $N$ sellers, and we shall drop the superscript $N$ from all the notations, for simplicity. 
Let $p^{*} = (p_1^{*}, p_2^{*}, \cdots, p_N^{*})$ be the candidate equilibrium prices (obtained by, 
for example,  the  first-order condition). By exchangeability, we can assume without loss of generality that the prices are ordered:
$p_1^{*} \leq p_2^{*} \leq \cdots \leq p_N^{*}$, and so are the corresponding cost functions 
$c_1^{*} \leq c_2^{*} \leq \cdots \leq c_N^{*}$,
where $c_i^{*} = c_i({p}^{*}, Q)$ for $i=1,\cdots, N$. 

%
%

We first  compare $p_N^{*,N}$ and $c_N^{*,N}$. 

{\bf Case 1.}  $p_N^{*} > c_N^{*}$. We consider the following cases:

(a) If $h_N^N (p^{*}) > 0$,  then by Definition \ref{ActualDemand} we have
$\h h_i(p^{*})=h^N_i(p^{*})>0$, for all $i$, and $ p^* = (p_1^{*}, p_2^{*}, \cdots, p_N^{*})$ is an
equilibrium point. 


(b) If $h_N^N (p^{*}) \le 0$, then in light of the definition of actual demand function (Definition \ref{ActualDemand}), we have
$\h h_N(p^{*})=0$. Thus, the $N$-th seller will have zero profit regardless where she sets the price. 
We shall require in this case that the $N$-th seller reduces her price to $c_N^{*}$, and we shall consider remaining $(N-1)$-sellers' 
candidate equilibrium prices $p^{*,N-1} = (p_1^{*,N-1}, \cdots, p_{N-1}^{*,N-1})$. 

{\bf Case 2.} $p_N^{*} \le c_N^{* }$. In this case the $N$-th seller would have a non-positive profit at the best. Thus, she sets $p_N^{*}= c_N^{* }$, and quits the game, and again the problem is reduced to a subgame with $(N-1)$ sellers, and in Case 1-(b). 

We should note that  it is possible that $h_N^N(p^*)> 0$ but $p_N^{*} \le c_N^{* }$. In this case, the profit would be non-positive, and the 
best option for the $N$-th seller is still to set  $p^*_N=c^*_N$ and quit. Such a case is known as  the ``boundary case", and we use
the notation  $p^*_N=c^{*,b}_N$ to indicate this case. 

Repeating the same procedure for the subgames (for $n=N-1, \cds, 2$), we see that eventually we will get a modified equilibrium point 
$p^*$ of the form (\ref{equilip}), proving the proposition. 
\qed

\subsection{A linear mean-field case}

In this subsection, we consider a special case, studied  in \cite{LS11}, but with the modified waiting cost functions. 
More precisely, we assume that there are $N$ sellers, each with demand function
\bea
\label{linearh}
h_i^N (p_1, \cdots, p_N) \triangleq A - B p_i + C \bar{p}_i^N ,
\eea
where $A,B,C >0$, and $B > C$, and $\bar{p}_i^N = \frac{1}{N-1} \sum\limits_{j \neq i}^{N} p_j$. We note that the structure of the 
demand function (\ref{linearh}) obviously reflects a mean-field nature, and one can easily check that it satisfies all the assumptions mentioned in the previous section. Furthermore, as was shown in \cite[Proposition 2.4]{LS11}, the actual demand function takes the form:
for each $n \in \{ 1, \cdots, N-1 \}$, 
\beaa
\label{actdemand}
h_i^n (p_1, \cdots, p_n) = a_n - b_n p_i + {c}_n \bar{p}_i^n, 
\q \text{ for } i=1, \cdots, n,
\eeaa
where 
$\bar{p}_i^n = \frac{1}{n-1} \sum\limits_{j \neq i}^{n} p_j$, and the parameters $(a_n, b_n, c_n)$ can be calculated recursively
for $ n=N, \cds, 1$, with  
$a_N = A$, $b_N = B$  and $c_N = C$. We note that in these works the (waiting) costs are assumed to be constant.

%
%

Let us now assume further that the waiting cost is also linear. For example,  for $n=1, \cdots, N$, 
\beaa
c_i^n = c_i^n (p_i, \bar{p}_i^n, Q) \triangleq x_n (Q) p_i - y_n (Q) \bar{p}_i^n, \qq x_n(Q), y_n(Q) >0. 
\eeaa

Note that the profit function for seller $i$ is 
\bea
\label{Pii}
\Pi_i (p_1, \cdots, p_n, Q) = (a_n - b_n p_i + c_n \bar{p}_i^n ) \cdot \left( p_i - (x_n p_i - y_n \bar{p}_i^n) \right). 
\eea
An easy calculation shows that the critical point for the maximizer is 
\bea
\label{pi*}
p_i^{*, n}= \frac{a_n}{2 b_n} + \left( \frac{{c}_n}{2 b_n} - \frac{y_n}{2(1-x_n)} \right) \bar{p}_i^n, 
\eea
 which is the optimal choice of seller $i$ if the other sellers set prices with average $\bar{p}_i^n = \frac{1}{n-1} \sum_{j\neq i}^{n} p_j^{*}$.
Now, let us define 
\bea
\label{barp}
\bar{p}^n := \frac{1}{n} \sum_{i=1}^{n} p_i^{*}=\frac{a_n (1-x_n)}{2b_n (1-x_n) - {c}_n (1-x_n) + b_n y_n}.
\eea
Then, it is readily seen that
$\bar{p}_i^n  = \frac{n}{n-1} \bar{p}^n - \frac{1}{n-1} p_i^*$, which means (plugging back into (\ref{pi*}))
\bea
\label{p*n}
p_i^{*, n} = \frac{a_n}{2b_n + \frac{{c}_n}{n-1} - \frac{1}{n-1} \frac{b_n y_n}{1-x_n}}
 + \frac{1}{\frac{n-1}{n} \frac{2b_n (1-x_n)}{{c}_n (1-x_n) - b_n y_n} + \frac{1}{n}} \bar{p}^n.
\eea

For the sake of argument, let us assume that the coefficients $(a_n, b_n, c_n, x_n(Q), y_n(Q))$ converge to $ (a,b,c,x(Q),
y(Q))$ as $n \rightarrow \infty$. Then, we see from (\ref{barp}) and (\ref{p*n}) that 
\bea
\label{limitp}
\left\{\ba{lll}
\dis \lim_{n\to\infty}\bar{p}^n =  \frac{a (1-x)}{2b(1-x)-c(1-x)+by}=:\bar{p};\ms\\
\dis \lim_{n\to\infty} p_i^{*, n} =\frac{a}{2 b} +  \frac{c(1-x)-by}{2b(1-x)} \lim_{n\to\infty}\bar{p}^n = \frac{a (1-x)}{(2b - c) (1-x) + b y}=:p^*.
\ea
\right.
\eea


 It is worth noting that if we 
 assume that there is a ``representing seller" who
 randomly sets prices $p=p_i$ with equal probability $\frac1n$, then we can randomize the profit function (\ref{Pii}):
\bea
\label{Pin}
\Pi_n(p, \bar p)= ( a_n - b_np + c_n \cdot \bar{p}) \{ p - (x_np - y_n\cd \bar{p} ) \},
\eea
where $p$ is a random variable taking values $\{p_i\}$ with equal probability, and $\bar{p}\sim \hE[p]$, thanks to the Law of Large Numbers, when $n$ is large enough. 
In particular, in the limiting case as $n\to \infty$,  we can replace the randomized profit function $\Pi_n$ in (\ref{Pin}) by:
\bea
\label{Pi}
\Pi_\infty=\Pi(p, \hE[p]):= ( a- bp + c  \hE[p]) \{ p - (xp - y \hE[p] ) \}.
\eea
A similar calculation as (\ref{pi*}) shows that $(p^*, \hE[p^*])\in \mbox{argmax} \Pi(p, \hE[p])$ will take the form
\begin{equation*}
p^* = \frac{c(1-x) - by}{2b (1-x)} \hE[p^*] + \frac{a}{2b}
\hspace{0.5cm} \text{ and } \hspace{0.5cm}
\hE[p^*]= \frac{a(1-x)}{2b(1-x) - c(1-x) + by}.
\end{equation*}
Consequently, we see that $p^* = \frac{a (1-x)}{(2b - c) (1-x) + b y}$, as we see in (\ref{limitp}). 

\begin{rem}
\label{remark2}{\rm
The analysis above indicates the following facts: (i) If we consider the sellers in a ``homogeneous" way, and as the 
number of sellers becomes large enough, all of them will actually choose the same strategy, as if there is a ``representing seller"
that places the prices uniformly; (ii) The limit of equilibrium prices actually coincides with the optimal strategy of the representing
seller under a limiting profit function. These facts are quite standard in mean-field theory, and will be used as the basis for our dynamic model for the (sell) LOB in the next
section.
\qed}
\end{rem}
%
%
%

\section{Mean-field type liquidity dynamics in continuous time}


In this section we extend the idea of Bertrand game to the  continuous time setting. To begin with, we assume that the contribution 
of each individual seller to the LOB is measured by the ``liquidity" (i.e., the number of shares of the given asset)
she provides, which is the function of the selling price she chooses, hence under the ``Bertrand game" framework. 

\subsection{A general description}
We begin by assuming that there are $N$ sellers, and denote the liquidity that the $i$-th seller ``adds" to the LOB at time $t$ by $Q_t^i$. We shall assume that it is a pure jump Markov process, with the following generator: for any $f\in C([0,T]\times \hR^N)$,  and $(t, q)\in [0,T]\times \hR^N$,
\bea
\label{generator-i}
\sA^{i,N}[f](t, q):=\int_{\hR}\l^i(t,q, \th)[f(t, q_{-i}(q_i+h^i(t,\th, z))-f(t,q)-\lan \pa_{x_i} f, h^i(t,\th, z)\ran]\n^i(dz),
\eea
where  $q\in \hR^N$, and $q_{-i}(y)=(q_1, \cds, q_{i-1}, y, q_{i+1}, \cds, q_N)$. Furthermore, 
$h^i$ denotes the ``demand function" for the $i$-th seller, and $\th\in \hR^k$ is a certain market parameter which will be specified
later. Roughly speaking,  (\ref{generator-i}) indicates that the $i$-th seller would act (or ``jump") at stopping times $\{\t^i_j\}_{j=1}^\infty$ with
the waiting times $\t^i_{j+1}-\t^i_j$ having exponential distribution with intensity $\l^i(\cd)$, and jump size being determined by the demand function $h^i(\cds)$. The total liquidity provided by all the sellers is then a pure jump process with the generator
\bea
\label{generator}
\sA^N[f](t,q, \th)=\sum_{i=1}^N \sA^{i, N}[f](t, q), \qq q\in\hR^N, ~N\in \hN,\qq  t\in [0,T].
\eea

We now specify the functions $\l^i$ and $h^i$ further. Recalling the  demand function introduced in the previous section, we assume that  there are
two functions $  \l$ and $ h$, such that for each $i$,
\bea
\label{tildelambdah}
\l^i(t, q, \th)=  \l(t, p^i,q^i, \m^N), \q h^i(t, \th, z)=  h(t, x, p^i, q^i,z), \qq (t,x,p, q)\in [0,T]\times \hR\times \hR^{2N}, 
\eea
where  $\m^N:=\frac1N\sum_{i=1}^N \d_{p^i}$, $x$ denotes the fundamental price at time $t$, and $p^i$ is the sell price. 
We shall consider $p=(p^1, \cds, p^N)$ as the control variable, as the Bertrand game suggests.  Now, if we assume  $\n^i=\n$ for all $i$, then we have
a {\it pure jump Markov game of mean-field-type}, similar to the one considered in \cite{BHK}, in which each seller adds liquidity
(in terms of number of shares) dynamically as a pure jump Markov process, denoted by $Q^i_t$, $t\ge0$, with the kernel
\bea
\label{kernel}
\n(t, q^i, \m^N, p^i, dz)=\l(t,  p^i, q^i,\m^N) [\n\circ h^{-1}(t, x, p^i,q^i, \cd )](dz).
\eea
Furthermore, in light of the static case studied in the previous section,  we shall assume that the seller's instantaneous profit at time $t>0$ takes the 
form $(p^i_t-c^i_t)\D Q^i_t$, where $c^i_t$ is the ``waiting cost" for $i$-th seller at time $t$. We observe that the actual 
submitted sell price $p^i$ can be written as $p^i=x+l^i$, where $x$ is the ``mid price" and $l^i$ is the distance from the mid price that the $i$-th seller chooses to set. Now let us assume that there is an invertible relationship between the selling prices $p$ and the corresponding number of shares $q$, e.g., $p =\f(q )$ (such a relation is often used to convert the Bertrand game to Cournot game, see, e.g., \cite{LS11}), and 
consider  $l$ as the {\it control variable}. We can then rewrite the functions $\l$ and $h$ of (\ref{tildelambdah}) in the following form: 
\bea
\label{lambdah}
\l(t, p^i,q^i, \m^N)= \l(t, q^i, l^i, \tilde\m^N(\f(q))), \q h(t,x, p^i, q^i,z)=h(t, x, q^i, l^i, z). 
\eea 
To simplify the presentation, in what follows, we shall assume that $\l$ does not depend on the control variable $l^i$, and  
that all coefficients $\lambda$ and $h$ are time-homogeneous. In other words, we assume that each $Q^i$ follows a pure jump SDE studied in \S2:
\bea
\label{Qi}
Q^i_t=q^i+\int_0^t\int_{A\times \hR_+} h(X_r, Q^i_{r-}, l^i_r,z){\bf 1}_{[0,\l(Q^i_{r-}, \m^N_{\f({\bf Q}_r)} )]}(y) \cN^s(drdzdy),  
\eea
where ${\bf Q}_t=(Q^1_t,\cds, Q^N_t)$,  $\cN^s$ is a Poission random measure on $\hR_+\times\hR\times\hR_+$, 
and $\{X_t\}_{t\ge 0}$ is the mid-price process of the underlying asset which we assume to  satisfy  the SDE (cf. \cite{MWZ}):
\begin{equation}
\label{midprice}
X_{s}^{t,x} = x + \int_{t}^{s} b(X_{r}^{t,x})dr + \int_{t}^{s} \sigma(X_{r}^{t,x})dW_r,
\end{equation}
where $b$ and $\sigma$ 
are deterministic functions satisfying some standard conditions. 
We shall assume that the $i$-th seller is aiming at maximizing the expected total accumulated profit:
\bea
\label{cost-i}
\hE\Big\{\sum_{t\ge 0} (p^i_t-c^i_t)\D Q^i_t\Big\}=\hE\Big\{\int_0^\infty\neg\neg\int_{A}\neg h(X_t, Q^i_t, l^i_t, z)(
X_t+l^i_t-c_t^i)\l(Q^i_t,\m^N_{\f({\bf Q}_t)}) \nu^s(dz)dt\Big\}.
\eea

We remark that in (\ref{cost-i}) the time horizon is allowed to be infinity, which can be easily converted to finite horizon by setting $h(X_t, \cds)=0$ for $t\ge T$, for a given time horizon $T>0$, which we do not want to specify at this point. Instead, our focus will be mainly on the limiting behavior
of the equilibrium when  $N\to \infty$. In fact, given the ``symmetric" nature of the problem (i.e., all seller's having the same $\l$ and
$h$), as well as the results in the previous section, we envision a ``representing seller"  in a limiting mean-field type control problem whose 
optimal strategy coincides with the limit of $N$-seller Nash equilibrium as $N\to \infty$, just as the well-known continuous diffusion cases (see, e.g., 
\cite{MFG} and \cite{BLPR, Carmona2}). We should note such a result for pure jump cases has been substantiated
in a recent work \cite{BHK}, in which it was shown that, under reasonable conditions, in the limit the total  liquidity $Q_t=\sum_{i=1}^N Q^i_t$ will converge to a pure jump Markovian process with a mean-field type generator. Based on this result, as well as the individual optimization problem
(\ref{Qi}) and (\ref{cost-i}), it is reasonable to consider the 
following (limiting) mean-filed-type pure-jump stochastic control problem for a representing seller, whose total liquidity has 
a dynamics that can be characterized by the following mean-field type pure-jump SDE:
\bea
\label{Q}
Q_t=q+\int_0^t\int_{A\times \hR_+} h(X_r, Q _{r-}, l _r,z){\bf 1}_{[0,\l(Q _{r-}, \hP_{Q_{r}})]}(y) \cN^s(drdzdy), 
\eea
where $\l(Q, \hP_Q):=\l(Q, \hE[\f(Q)])$ by a slight abuse of notation,  and with the cost functional:
\bea
\label{cost}
\Pi(q, l)=\hE\Big\{\int_0^\infty\neg\neg\int_{A}\neg h(X_t, Q _t, l_t, z)(
X_t+l_t-c_t)\l(Q_t,\hP_{{Q}_t}) \nu^s(dz)dt\Big\}.
\eea

\subsection{Problem formulation}

With the general description in mind, we now give the formulation of our problem. First, we note that the liquidity of the limit order book will
not only be affected by the liquidity providers (i.e., the sellers), but also by liquidity ``consumer", that is, the market buy orders as well as 
the cancellations of sell orders (which we assume is free of charge). We shall describe its collective movement (in terms of number of shares) 
of all such ``consumptional" orders as a compound Poisson process, denoted by $\b_t=\sum_{i=1}^{N_t} \L_t$, $t\ge 0$, where $\{N_t\}$ is 
a standard Poisson process with parameter $\l$, and $\{\L_i\}$ is a sequence of i.i.d. random variables taking values in a set
$B\subseteq \hR$, with distribution $\n$. Without loss of generality, we assume that counting measure of $\b$ coincides with the canonical 
Poisson random measure $\cN^b$, so that the L\'evy measure is $\n^b=\l \n$. In other words, $\b_t:= \int_0^t \int_B z \; \wt\cN^b(dr dz)$, and the total liquidity satisfies the SDE:
\bea
\label{Qb}
\label{Q}
Q^0_t=q+\int_0^t\int_{A\times \hR_+} h(X_r, Q^0_{r-}, l _r,z){\bf 1}_{[0,\l(Q^0_{r-}, \hP_{Q^0_{r}})]}(y) \cN^s(drdzdy)-\b_t. 
\eea

We remark that there are two technical issues for the dynamics (\ref{Qb}). First, the presence of the buy order process $\b$ brings in the possibility
that $Q^0_t<0$, which should never happen in reality. We shall therefore assume that  the buy order has a natural upper limit: the total available liquidity $Q^0_t$, that is, if we denote $\cS_\b=\{t:\D\b_t\neq 0\}$, then for all $t\in \cS_\b$, we have $Q^0_t=(Q^0_{t-}-\D \b_t)^+$. Consequently, we can assume that there exists a process $K=\{K_t\}$, where $K$ is a non-decreasing, pure jump process such that (i) $\cS_K=\cS_\b$; (ii) $\D K_t:=(Q_{t-}-\D \b_t)^-$, $t\in \cS_K$; and (iii) the $Q^0$-dynamics (\ref{Qb}) can be written as 
\bea
\label{QbK}
Q_t&=&q+\int_0^t\neg\neg\int_{A\times \hR_+} h(X_r, Q _{r-}, l _r,z){\bf 1}_{[0,\l(Q _{r-}, \hP_{Q_{r}})]}(y) \cN^s(drdzdy)-\b_t +K_t\nonumber\\
&=& q+\int_0^t\int_{A\times \hR_+} h(X_r, Q _{r-}, l _r,z){\bf 1}_{[0,\l(Q _{r-}, \hP_{Q_{r}})]}(y) \wt \cN^s(drdzdy) \\
&&-\int_0^t\neg\neg \int_B z \; \wt\cN^b(dr dz)+\int_0^t \neg\neg\int_{A} h(X_r, Q _{r}, l _r,z)\l(Q _{r}, \hP_{Q_{r}})\n^s(dz)dr+K_t, \q t\ge 0. \nonumber
\eea
where $K$ is a ``reflecting process", and $\wt\cN^s(drdzdy)$ is the compensated Poisson martingale measure of $\cN^s$. That is,  
%
%
%
(\ref{QbK})  is a (pure-jump) mean-field SDE with reflection as was studied in \S2.  

Now, in light of the discussion of MFSDER in \S2, we shall consider the following two MFSDERs that are slightly more general than (\ref{QbK}): 
for $\xi \in L^2(\cF_t; \dbR)$,  $q \in \dbR$, and $0\le s\le t$, 
\bea
\label{QL2}
Q_{s}^{t,\xi} &=& \xi 
+ \int_t^s \int_{A \times \dbR_{+}} h(X_r^{t,x}, Q_{r-}^{t,\xi},  l_r, z)  \;
\1_{ [0, \lambda(Q_{r-}^{t,\xi}, \hP_{Q_r^{t,\xi}})] } (y) \wt\cN^{s}(dr dz dy) \nonumber \\
&&-  \int_t^s \int_B z \; \wt\cN^b(dr dz)
+ \int_t^s a(X_r^{t,x}, Q_r^{t,\xi}, \hP_{Q_{r}^{t,\xi}}, l_r) dr
+ K_s^{t,\xi}, \\
\label{Qq}
Q_{s}^{t,q,\xi} &=& q 
+ \int_t^s \int_{A \times \dbR_{+}} h(X_r^{t,x}, Q_{r-}^{t,q,\xi},  l_r, z) \;
\1_{ [0, \lambda(Q_{r-}^{t,q,\xi}, \hP_{Q_r^{t,\xi}})] } (y) \wt\cN^{s}(drdzdy) \nonumber \\
&& - \int_t^s \int_B z \; \wt\cN^b(dr dz) 
+ \int_t^s a(X_r^{t,x}, Q_r^{t,q,\xi}, \hP_{Q_{r}^{t,\xi}}, l_r) dr 
+ K_s^{t,q,\xi},
\eea
where $l=\{l_s\}$ is the control process for seller, and $Q_s=Q^{t, q, \xi}_s$, $s\ge t$, is the total liquidity of the sell-side LOB.
We shall consider the following set of  \textit{admissible strategies}:
\bea
\label{Adm}
\sU_{ad}:=\{ l \in L^1_\hF([0,\infty);\hR_{+}): l \text{ is } \dbF \text{-predictable}\}.
\eea

The objective of the seller is to solve the following mean-field stochastic control problem:
\bea
\label{value}
 v(x,q,\hP_\xi) =  \sup_{l \in \sU_{ad}} \Pi (x,q,\hP_\xi, l) = \sup_{l \in \sU_{ad}} \dbE \Big\{
\int_0^\infty  e^{-\rho r} L(X_r^{x}, Q_r^{q,\xi}, \hP_{Q_r^{\xi}}, l_r)  dr
 \Big\}
\eea
where $L(x, q, \m, l)  := \int_A h(x, q, l, z) c(x, q, l) \l(q, \m)\nu^s (dz)$, and $\sU_{ad}$ is defined in (\ref{Adm}). 
Here we denote $X^x := X^{0,x}$, $Q^{q, \xi} := Q^{0,q,\xi}$. 

\begin{rem}
{\rm
(i) In (\ref{QL2}) and (\ref{Qq}), we allow a slightly more general drift function $a$, which in particular could be  $a(x, q, \m, l) = \lambda(q, \m) \int_A h(x, q, l, z) \nu^s (dz)$, as is in (\ref{QbK}). 


\ms
(ii) In (\ref{value}), the pricing function $c(x,q,l)$ is a more general expression of  the original form $x+l-c$ in (\ref{cost}), taking into account the possible dependence of the waiting cost $c_t$  on the sell position $l$   and
the total liquidity $q$ at time $t$. 

\ms
(iii) Compared to (\ref{cost}), we see that a {\it discounting factor} $e^{-\rho t}$ is added to the cost functional $\Pi(\cds)$ in (\ref{value}), reflecting its nature as the ``present value".
\qed
}
\end{rem}

In the rest of the paper we shall  assume that the market parameters $b, \si, \l, h$, the pricing function $c$ in (\ref{QL2}) -- (\ref{value}), and the discounting factor $\rho$ satisfy 
the following assumptions. 

\begin{assum}
\label{assump2}
All functions $b, \si \in C^0(\hR)$, $\l \in L^0(\mathbb{R} \times \sP_{2}(\mathbb{R}); \mathbb{R}_+)$,  
$h\in L^0(\hR^2\times \hR_+ \times A)$, and $c\in L^0( \hR\times\hR_+ \times \hR_+)$ are bounded, 
and satisfy the following conditions, respectively:

%

\ss
(i) $b$ and $\si$ are  uniformly Lipschitz continuous in $x$  with Lipschitz constant $L>0$;

\ss
(ii) $\sigma(0)=0$ and  $b(0) \geq 0$; 

\ss
(iii) $\l$ and $h$ satisfy Assumption \ref{assump1};


\ss 
(iv) For $l \in \hR_+$, $c(x,q,l)$ is Lipschitz continuous in $(x,q)$, with Lipschitz constant $L>0$;

\ss
(v) $h$ is non-increasing, and $c$ is non-decreasing in the variable $l$;

\ss
(vi)  $\rho> L+\frac12 L^2$, where $L>0$ is the Lipschitz constant in Assumption \ref{assump1};

\ss
(vii) For $(x,\mu, l) \in \dbR_{+} \times \sP_2 (\dbR) \times \dbR_+$, $\Pi(x,q,\mu, l)$ is convex in $q$. 
 \qed
 \end{assum}



\begin{rem}
\label{remark3}
{
\rm (i) The monotonicity assumptions in Assumption \ref{assump2}-(v) are inherited from 
\S3. Specifically, they are the assumption (\ref{mono})  for $h$, and Assumption 3.1-(i) for $c$, respectively. 

(ii) Under Assumption 4.2, one can easily check that the SDEs (\ref{midprice}) as well as (\ref{QL2}) and (\ref{Qq}) all have pathwisely unique strong solutions in $L^2_\hF (\hD([0,T]))$, thanks to Theorem \ref{thm1-MFSDER}; and Assumption \ref{assump2}-(ii) implies that $X^{t, x}_s\ge 0$, $s\in[t, \infty)$, $\hP$-a.s., whenever $x\ge 0$.
\qed}
\end{rem}

%
%
%
%


\section{Dynamic programming principle and HJB equation}

In this section we substantiate the {\it dynamic programming principle} (DPP) for the stochastic
control problem (\ref{QL2})--(\ref{value}).
We begin by examining some basic properties of the value function. 
\begin{proposition}
\label{Prop_ValueFunc}
Under the Assumptions \ref{assump1} and \ref{assump2}, 
the value function $v(x,q, \hP_{\xi})$ is 
Lipschitz continuous in $(x,q,\hP_{\xi})$, 
{non-decreasing in $x$, and decreasing in $q$}.
\end{proposition}

{\it Proof.} 
We first check the Lipschitz property in $x$. 
For $x, x'\in\hR$, denote $X^x=X^{0,x}$ and $X^{x'}=X^{0,x'}$ as the corresponding solutions to (\ref{midprice}), respectively. Denote $\D X_t=X^x_t-X^{x'}_t$, and $\D x=x-x'$. Then,  applying
 It\^o's formula to $|\D X_t|^2$ and by some standard arguments, one has 
$$ |\D X_t|^2=|\D x|^2+\int_0^t (2\a_s +\b_s^2) |\D X_s|^2ds+\int_0^t 2\b_s |\D X_s|^2dW_s,
$$
where $\a, \b$ are two processes bounded by the Lipschitz constants $L$ in Assumption \ref{assump1}, thanks to Assumption \ref{assump2}. Thus, one can easily check, by taking expectation and applying Burkholder-Davis-Gundy and Gronwall inequalities, that 
\bea
\label{Xest}
\hE[|\D X|^{*,2}_t]\le |\D x|^2 e^{(2L+L^2)t}, \qq t\ge 0.
\eea

Furthermore, it is clear that, under Assumption \ref{assump2}, the function $L( x,q, \m, l)$ is uniformly 
Lipschitz in $x$, uniformly in $(q,\m, l)$. That is, for some generic constant $C>0$ 
which is allowed to vary from line to line, we have
\beaa
&&|\Pi (x,q,\hP_{\xi},l) - \Pi (x',q,\hP_{\xi},l)| \le C\hE \Big[
\int_0^\infty \int_A  e^{-\rho t}
|\D X_t |   \nu^s (dz) dt\Big]  \\
&\leq&  C  \hE \Big[\int_0^\infty e^{-\rho t} \sqrt{\hE[|\D X|^{*,2}_t]}dt\Big]\le C|\D x| \int_0^\infty e^{-\rho t}e^{(L+\frac12 L^2)t}dt
\leq C |x - x'|.
\eeaa
Here the last inequality is due to Assumption \ref{assump2}-(vi). Consequently, we obtain
\begin{equation}
 | v(x, q, \hP_{\xi}) - v(x',  q, \hP_{\xi}) | \leq C |x-x'|, 
 \hspace{0.3cm} \forall x,x' \in \dbR.
\end{equation}

To check the Lipschitz properties for $q$ and  $\hP_\xi$, we denote, for $(q, \hP_\xi)\in \dbR_{+}\times \sP_{2}(\dbR)$, 
$h_s^{q, \xi} \equiv h(X_s, Q_s^{t,q, \xi}, l_s, z)$, 
$\lambda_s^{q, \xi} \equiv \lambda (Q_s^{t,q, \xi}, \hP_{Q_s^{t,\xi}})$, and
$c_s^{q,\xi} \equiv c(X_s, Q_s^{t,q, \xi}, l_s)$, $s \ge t$. 
Furthermore, for  $q, q' \in \dbR_{+}$ and $\hP_\xi, \hP_{\xi'} \in \sP_{2}(\dbR)$,
we deonte $\Delta \psi_r \equiv  \psi_r^{q,\xi} - \psi_r^{q', \xi'}$
for $\psi = h, \lambda, c$.
Now, by Assumptions \ref{assump1} and \ref{assump2}, and following a similar argument of Theorem \ref{thm1-MFSDER}, one shows
that
\beaa
&&| \Pi(x, q, \hP_\xi, l) - \Pi (x,  q', \hP_{\xi'}, l) | \\
&\le &  \hE \Big\{\int_0^\infty\neg \int_A e^{-\rho r} \big(  h_r^{q,\xi} c_r^{q,\xi}|\Delta \lambda_r |
+ c_r^{q,\xi} \lambda_r^{q', \xi'} |\Delta h_r|+ h_r^{q',\xi'} \lambda_r^{q', \xi'} |\Delta c_r |  \big)  \nu^s (dz) dr
\Big\} \\
\neg&\neg \leq \neg&\neg \hE \Big\{\int_0^\infty\neg \int_A e^{-\rho r} |Q^{t,q,\xi}_r-Q_r^{t,q',\xi'} |  \nu^s (dz) dr
\Big\} \le C \left( |q-q'| + W_1(\hP_\xi,\hP_{\xi'}) \right),
\eeaa
which implies that 
\begin{equation}
 | v(x, q, \hP_\xi) - v(x,  q', \hP_\xi') | 
 \leq C \left( |q-q'| +W_1 (\hP_\xi, \hP_{\xi'}) \right).
\end{equation}

Finally, the respective monotonicity of the value function on $x$ and $q$   follows from the comparison theorem of the
corresponding SDEs  and Assumption
\ref{assump2}. This completes  the proof. 
\qed

We now turn our attention to the Dynamic Programming Principle (DPP). 
The argument will be very similar to that of \cite{MWZ}, except for some adjustments to deal with the mean-field terms. But, by using the flow-property (\ref{flow}) we can carry out the argument  without substantial difficulty.
\begin{thm}
\label{DPP}
Assume that Assumptions \ref{assump1} and \ref{assump2} are in force. 
Then, for any $(x, q, \hP_{\xi}) \in \dbR^2 \times
\sP_2 (\dbR)$ and for any $t \in (0, \infty)$, 
\bea
\label{DPP_eqn}
v(x,q,\hP_{\xi}) = \sup_{l \in \sU_{ad}} \dbE \Big[ \int_0^t e^{-\rho s} L(X_s^{x}, Q_s^{q,\xi; l}, \hP_{Q_{s}^{\xi; l}}, l_s) ds  
+ e^{-\rho t} v( X_{t}^{x}, Q_{t}^{q,\xi; l}, \hP_{Q_{t}^{\xi ;l}}) \Big]. 
\eea
\end{thm}

{\it Proof.}
Let us denote the right side of (\ref{DPP_eqn}) by $\tilde{v}(x,q, \hP_{\xi}) = \sup_{l} \tilde{\Pi} (x,q, \hP_{\xi};l)$. 
We first note that $X_r$ and $(Q_r^{t,\xi}, Q_{r}^{t,q,\xi})$ have the flow property. So, for any $l\in \sU_{ad}$, 
\bea
\label{DPPest1}
 &&\Pi(x,q,\hP_\xi; l) =  \hE \Big[ \int_0^\infty e^{-\rho s} L(X_s^{x}, Q_s^{q,\xi; l}, \hP_{Q_{s}^{\xi; l}}, l_s) ds \Big] \\
\neg \neg&\neg = \neg&\neg\neg \hE \Big[ \int_{0}^{t} e^{-\rho s} L(X_s^{x}, Q_s^{q,\xi; l}, \hP_{Q_{s}^{\xi; l}}, l_s) ds
+ e^{-\rho t} \hE \Big\lbrace \int_{t}^{\infty} e^{-\rho (s- t)} L( X_s^{x}, Q_s^{q,\xi ;l}, \hP_{Q_s^{\xi;l}}, l_s) ds
\Big|\cF_{t} \Big\rbrace\Big] \nonumber \\
\neg \neg&\neg = \neg&\neg \neg \hE \Big[ \int_{0}^{t} e^{-\rho s} L(X_s^{x}, Q_s^{q,\xi; l}, \hP_{Q_{s}^{\xi; l}}, l_s) ds 
+ e^{-\rho t}  \Pi(X_t^{x}, Q_t^{q,\xi ;l}, \hP_{Q_t^{\xi;l}}; l) \Big] \nonumber\\
\neg\neg&\neg \leq \neg&\neg\neg  \hE \Big[ \int_{0}^{t} e^{-\rho s} L(X_s^{x}, Q_s^{q,\xi; l}, \hP_{Q_{s}^{\xi; l}}, l_s) ds 
+ e^{-\rho t} v(X_t^{x}, Q_t^{q,\xi ;l}, \hP_{Q_t^{\xi;l}}) \Big] \nonumber \\
\neg\neg&\neg = \neg&\neg\neg \tilde{\Pi} (x,q,\hP_{\xi};l). \nonumber
\eea
This implies that  $v( x,q,\hP_{\xi}) \leq \tilde{v}(x,q,\hP_{\xi})$. 

To prove the other direction, 
let us denote $\Gamma = \hR_{+} \times \hR \times \sP_2 (\dbR)$, and consider,  
at each time $t\in(0,\infty)$, a countable partition $\{ \Gamma_i \}_{i=1}^{\infty}$ of $\Gamma$ and $(x_i, q_i, \hP_{\xi_i}) \in \G_i$, $\xi_i\in L^2(\cF_{t})$, $i=1, 2, \cds$,  such that for any $(x,q,\mu) \in \G_i$ and for fixed $\e >0$, it holds
$ |x-x_i | \leq \e $, $q_i - \e \leq q \leq q_i$, and $W_2 (\mu, \hP_{\xi_i}) \leq \e$.
Now, for each $i$, choose an $\e $-optimal strategy $l^i \in \sU_{ad}$, such that 
$v(t,x_i, q_i, \hP_{\xi_i}) \leq \Pi(t, x_i, q_i, \hP_{\xi_i}; l^i) + \e$, 
where 
$\Pi (t, x_i, q_i, \hP_{\xi_i}; l^i) := 
\hE [\int_t^{\infty} e^{-\rho (s-t)} L(X_s^{t, x_i}, Q_s^{t, q_i, \xi_i}, \hP_{Q_{s}^{t, \xi_i}}, l^i_s) ds]$ 
and $v(t,x_i, q_i, \hP_{\xi_i}) = \sup_{l^i \in \sU_{ad}} \Pi (t, x_i, q_i, \hP_{\xi_i}; l^i)$.

Then, by definition of the value function and the Lipschitz properties (Proposition \ref{Prop_ValueFunc}) with 
some constant $C>0$,  for any $(x, q, \m)\in \G_i$, it holds that 
\bea
\label{DPPest2}
\Pi(t, x,q,\mu; l^i) &\geq& \Pi(t, x_i, q_i, \hP_{\xi_i}; l^i) - C \e 
\geq v(t, x_i, q_i, \hP_{\xi_i}) -( C+1) \e \nonumber\\
&\geq& v(t, x,q,\mu) -(2C+1)\e.  
\eea

Now, for any $l \in \sU_{ad}$, we define a new strategy $\tilde{l}$ as follows: 
\bea
\tilde{l}_s := l_s \1_{ [0,t] } (s) 
+ \Big[ \sum_{i} l^i_s \1_{ \G_i } (X_{t}^{ x}, Q_{t}^{ q, \xi; l}, \hP_{Q_{t}^{\xi; l}}) \Big]
\1_{ (t, \infty) } (s).
\eea
Then, clearly $\tilde{l} \in \sU_{ad}$. 
To simplify notation, let us denote
\bea
\label{I1}
I_1 = \int_{0}^{t} e^{-\rho s} L(X_s^{x}, Q_s^{q,\xi; l}, \hP_{Q_{s}^{\xi; l}}, l_s) ds.
\eea
By applying (\ref{DPPest2}) and flow property, we have
\beaa
&&v(x,q,\mu) \ge   \Pi(x,q,\mu; \tilde{l}) \\
\neg \neg&\neg = \neg&\neg\neg \dbE \Big[ I_1 
+ e^{-\rho t}   \hE \Big\lbrace \int_{t}^{\infty} e^{-\rho (s- t)} L(X_s^{x}, Q_s^{q,\xi; l}, \hP_{Q_{s}^{\xi; l}}, l_s)ds
\Big|\cF_{t} \Big\rbrace \Big]  \\
\neg \neg&\neg = \neg&\neg \neg \dbE \Big[ I_1
+ e^{-\rho t}  \Pi(t,X_{t}^{x}, Q_{t}^{q,\xi}, \hP_{Q_{t}^{\xi}}; \tilde{l}
) \Big] \\
\neg\neg&\neg = \neg&\neg\neg  \dbE \Big[ I_1
+ e^{-\rho t} \sum_i \Pi(t, X_{t}^{x}, Q_{t}^{q,\xi}, \hP_{Q_{t}^{\xi}}; l^i) 
\1_{\Gamma_i} (X_{t}^{x}, Q_{t}^{q,\xi},\hP_{Q_{t}^{\xi}}) \Big] \\
\neg\neg&\neg \geq \neg&\neg\neg \dbE \Big[ I_1
+ e^{-\rho t} v(X_{t}^{x}, Q_{t}^{q,\xi},\hP_{Q_{t}^{\xi}}) \Big]  -(2 C+1) \e = \tilde{\Pi} (x,q,\hP_{\xi};l) -(2 C+1) \e.
\eeaa
\no Since $\e >0$ is arbitrary, we get $v(x,q,\hP_{\xi}) \geq \tilde{v}(x,q,\hP_{\xi})$, proving
(\ref{DPP_eqn}). 
\qed

{
\begin{rem}
\label{q=0}
{\rm
We should note that while it is difficult to specify all the boundary conditions for the value function, the case when $q=0$ is relatively clear. 
Note that $q=0$ means there is zero liquidity for the asset, then by definition of the liquidity dynamics (\ref{QbK}) we see that $Q_t$ will 
stay at zero until the first positive jump happens. During that period of time there would be no trade, thus by DPP (\ref{DPP_eqn}) we should have
 \bea
\label{bdyv}
v(x,0, \m)\equiv 0.
\eea
Furthermore, since the value function $v$ is non-increasing in $q$, thanks to Proposition \ref{Prop_ValueFunc}, and is always
non-negative, we can easily see that the following boundary condition is also natural. \bea
\label{bdyvq}
\pa_q v(x, 0,\m)\equiv 0.
\eea
We shall use (\ref{bdyv}) and (\ref{bdyvq}) frequently in our future discussion.
 \qed}
\end{rem}
}

\section{HJB equation and its viscosity solutions}

In this section, we shall formally derive the Hamilton-Jacobi-Bellman (HJB) equation associated to the stochastic control problem studied in the 
previous section, and show that the value function of the control problem is indeed a viscosity solution of the HJB equation. 

To begin with,  we first note that, given the DPP (\ref{DPP_eqn}), {
as well as the boundary conditions (\ref{bdyv}) and (\ref{bdyvq}), if the value function $v$ is smooth, then by standard arguments with the help
of the It\^o's formula (\ref{ItoForm}) and the fact that 
$$\pa_q v(X_{t-}, Q_{t-}, \hP_{X_{t-}}) \1_{\{Q_{t-}=0\}}dK_t =\pa_q v(X_{t-}, 0, \hP_{X_{t-}}) \1_{\{Q_{t-}=0\}}dK_t\equiv 0,
$$
 it is not difficult to show that the value function should satisfy the following HJB equation: 
\bea
\label{HJB_seller}
\left\{\ba{lll}
\dis \rho v (x,q,\m) = \sup_{l \in\hR_+} [ \sJ^{l} v(x,q,\m) + L(x,q,\m,l) ], \ms\\
%
v(x,0,\m) =0, \q {\pa_q v(x, 0, \m)=0,}
\ea\right.  \q (x,q, \m)\in \hR\times\hR_+\times \sP_2(\hR),
\eea
where $\sJ^l$ is an integro-differential operator defined by, for any $\phi \in \hC_b^{2,(2,1)} (\dbR \times\dbR_{+} \times  \sP_2(\dbR))$,  
\bea
&&\sJ^{l} \phi (x,q, \m)  \dfnn 
b(x) \partial_x \phi (x,q,\m) + \sigma^2 (x) \frac{1}{2} \partial_{xx}^{2} \phi(x,q,\m) + a(x,q, \mu, l) \partial_q \phi(x,q,\m)  
\nonumber\\
&&\q +\int_{A} \Big( \phi(x, q+ h(x, q,  l, z), \m ) - \phi(x,q,\m) - \partial_q \phi(x,q,\m) h(x,q,l,z)\Big)
\lambda(q,\m) \nu^{s}(dz)  \nonumber\\
&&\q - \int_B \Big(
\phi(x,q-z, \m) - \phi(x,q, \m) - \partial_q \phi(x,q, \m) z
\Big) \nu^{b} (dz) 
+ \tilde{\dbE} \left[ \partial_\mu \phi (x,q, \m, \tilde{\xi}) a(x, \tilde{\xi}, \m, l) \right]  \nonumber\\
&&\q + \tilde{\dbE} \Big[ \int_0^1\neg\neg \int_{A} \Big( \partial_\mu \phi (x,q, \m, \tilde{\xi}+ \gamma h(x,\tilde{\xi},  l, z) )
 - \partial_\mu \phi (x,q, \m, \tilde{\xi}) \Big)
 h(x,\tilde{\xi}, l,z) \lambda(\tilde{\xi}, \m) \nu^{s} (dz) d \gamma \Big]  \nonumber \\
&&\q - \tilde{\dbE} \Big[ \int_0^1\neg\neg \int_B \Big( \partial_\mu \phi(x,q, \m, \tilde{\xi}- \gamma z) - \partial_\mu \phi(x,q, \m, \tilde{\xi}) \Big) \times z 
\nu^{b}(dz) d\gamma \Big].
\eea
}
We note that in general, whether there exists smooth solutions to the HJB equation (\ref{HJB_seller}) is by no means clear. We therefore 
introduce the notion of {\it viscosity solution} for  (\ref{HJB_seller}). 
To this end, write $\sD :=\dbR \times \hR_{+} \times  \sP_2 (\dbR)$, and for $(x,q,\mu) \in \sD$, we denote
\beaa
\sU (x,q, \mu) &:= & \Big\{ \varphi \in \hC_b^{2,(2,1)} (\sD) : v(x,q,\mu) = \varphi (x,q,\mu) \Big\}; \nonumber \\
\overline{\sU}(x,q,\mu) &:= & \Big\{ \varphi \in \sU (x,q, \mu)  : v- \varphi \text{ has a strict maximum at } (x,q,\mu) \Big\}; \\
\underline{\sU}(x,q,\mu) &:= & \Big\{ \varphi \in \sU (x,q, \mu)  : v- \varphi \text{ has a strict minimum at } (x,q,\mu) \Big\}. \nonumber
\eeaa
\begin{defn}
\label{Def_viscosity}
We say a continuous function $v: \sD \mapsto \dbR_{+}$ is a \textit{viscosity subsolution (supersolution, resp.)} of (\ref{HJB_seller}) in $\sD$ if
\begin{equation}
\rho \varphi(x,q,\mu) - \sup\limits_{l \in \hR_+} [ \sJ^{l} \varphi(x,q,\mu) +  L(x,q,\mu,l)
 ] \leq 0, ( \mbox{resp.} \geq 0)
\end{equation}
for every $\varphi \in \overline{\sU}(x,q,\mu)$ (resp. $\varphi \in \underline{\sU}(x,q,\mu))$.

A function $v:\sD\mapsto \hR_+$ is called a viscosity solution of (\ref{HJB_seller}) on $\sD$ if it is both a viscosity subsolution and a viscosity supersolution of (\ref{HJB_seller}) on $\sD$.
 \qed
\end{defn}

%
%
%

Our main result of this section is  the following theorem.

\begin{theorem}
Assume that the Assumptions \ref{assump1} and \ref{assump2} are in force. Then, the value function $v$, defined by (\ref{value}), is a viscosity solution of the HJB equation (\ref{HJB_seller}).
\end{theorem}
{\it Proof.}
For a fixed  $\bar\dbx := (\bar x, \bar q,\bar\mu) \in \sD$ with $\bar\mu = \hP_{\bar{\xi}}$ and $\bar\xi\in L^2(\cF; \dbR)$,  and any $\eta >0$, define
\bea
\label{Dxeta}
\sD_{\bar\dbx, \eta}:=\{\dbx = (x, q, \mu)\in \sD: \|\dbx-\bar\dbx\|<\eta\}.
\eea
where  $\| \dbx - \bar{\dbx} \| := \Big( |x-\bar{x}|^2 +|q-\bar{q}|^2 +W_2 (\mu, \bar{\mu})\Big) ^{1/2}$, and $\mu = \hP_{\xi}$ with $\xi \in L^2(\cF; \dbR)$.  

We first prove that the value function $v$ is a {\it subsolution} to the HJB equation (\ref{HJB_seller}).
We proceed by contradiction. Suppose not. Then there exist some $\f\in \ol\sU(\dbx)$ and $\e_0>0$ such that
\begin{equation}
\rho \f(\bar\dbx) - \sup_{l \in \hR_+} [ \sJ^{l} [\f](\bar\dbx) + L(\bar\dbx,l) ]  =: 2 \e_0 > 0.
\end{equation}
Since $A^l(\bmx):= \sJ^{l} \varphi(\bmx) + L(\bmx,l)$ is uniformly continuous in $\bmx$, uniformly in $l$, thanks to Assumption \ref{assump2}, one
shows that there exists $\eta >0$ such that for any $\dbx \in \sD_{\bar\dbx, \eta}$, it holds  that 
\bea
\label{epsilon}
\rho \varphi(\dbx) - \sup_{l \in \hR_+} [ \sJ^{l} [\varphi](\dbx) + L(\dbx,l) ]  \geq \e_0 .
\eea
Furthermore, since $\f\in\ol\sU(\dbx)$, we assume without loss of generality that $0=v(\bar\bmx)-\f(\bar\bmx)$ is the strict maximum. Thus for the given $\eta>0$, there exists $\delta >0$, such that
\bea
\label{delta}
\max \left\lbrace v( {\dbx}) - \varphi( {\dbx}): {\dbx} \notin \sD_{\bar\dbx, \eta}\right\rbrace 
= -\delta < 0,
\eea

On the other hand, for a fixed $\e \in (0, \min(\e_0, \delta \rho))$, by the continuity of $v$ we can assume, modifying $\eta>0$ if necessary, that 
\bea
\label{veta}
|v(\dbx)-v(\bar\dbx)|=|v(\dbx)-\f(\bar \dbx)|<\e, \qq \dbx\in \sD_{\bar\dbx, \eta}.
\eea

Next, for any $T>0$  and any $l\in \sU_{ad}$ we set $ \tau^T := \inf \{ t \geq 0 : \bar\Theta_t \notin \sD_{\bar\dbx, \eta} \} 
\wedge T$, where  $\bar\Theta_t:= (X_{t}^{\bar x}, Q_t^{\bar q, \bar\xi, l}, \hP_{Q_t^{\bar q, \bar\xi}} )$.
%
Applying It\^o's formula (\ref{ItoForm}) to $e^{-\rho t} \varphi(\bar\Theta_t)$ from $0$ to $\t^T$ and noting that $v(\bar\dbx)=\f(\bar\dbx)$ we have 
\bea
\label{estbar}
&&\dbE \Big[ \int_0^{\t^T} e^{-\rho t} L(\bar\Theta_t, l_t) dt 
+e^{-\rho \t^T} v(\bar\Theta_{\t^T} )\Big]\nonumber\\
&=&\dbE \Big[ \int_0^{\t^T} e^{-\rho t} L(\bar\Theta_t, l_t) dt 
+e^{-\rho \t^T} \f(\bar\Theta_{\t^T} )+e^{-\rho \t^T} [v-\f](\bar\Theta_{\t^T} )\Big]\nonumber\\
&=& \dbE \Big[ \int_0^{\t^T} e^{-\rho t} \left( L(\bar\Th_t, l_t) + \sJ^{l} [\varphi](\bar\Theta_{t}) - \rho \varphi(\bar\Th_{t})
\right) dt +e^{-\rho \t^T} [v-\f](\bar\Theta_{\t^T} ) \Big]+v(\bar\dbx)  \\
&\le &\hE\Big[ - \frac{\e}{\rho}(1-e^{-\rho \t^T})+e^{-\rho \t^T} [v-\f](\bar\Theta_{\t^T} ) \Big]+v(\bar\dbx) \nonumber\\
&=&\hE\Big[ e^{-\rho \t^T} \Big(\frac{\e}{\rho}+[v-\f](\bar\Theta_{\t^T})\Big):\t^T<T \Big]+\hE\Big[e^{-\rho \t^T} \Big(\frac{\e}{\rho}+[v-\f](\bar\Theta_{\t^T})\Big):\t^T=T \Big] \nonumber\\
&&+v(\bar\dbx)- \frac{\e}{\rho}.\nonumber
\eea
Now note that on the set $\{\t^T<T\}$ we must have $\bar\Th_{\t^T}\notin \sD_{\bar\dbx, \eta}$, thus $[v-\f](\bar\Th_{\t^T})\le -\d$, thanks to 
(\ref{delta}). On the other hand, on the set $\{\t^T=T\}$ we have $\bar\Th_{\t^T}=\bar\Th_T\in \sD_{\bar\dbx,\eta}$, and then (\ref{veta}) implies that
$[v-\f](\bar\Theta_{T})\le v(\bar\dbx)-\f(\bar\Th_T)+\e$. Plugging these facts in (\ref{estbar}), we can easily obtain that
\beaa
\dbE \Big[ \int_0^{\t^T} e^{-\rho t} L(\bar\Theta_t, l_t) dt +e^{-\rho \t^T} v(\bar\Theta_{\t^T} )\Big] 
&\le& \Big(\frac{\e}{\rho}-\d\Big)\hP\{\t^T<T\}+(\frac{\e}{\rho}+\e)e^{-\rho T}+v(\bar\dbx)-\frac{\e}{\rho}\\
&\le&\Big(\frac{\e}{\rho}+\e\Big)e^{-\rho T}+v(\bar\dbx)-\frac{\e}{\rho}.
\eeaa
Here in the last inequality above we used the fact that $\e/\rho -\d<0$, by definition of $\e$. 
Letting $T \rightarrow \infty$ we have 
\beaa
\dbE \Big[ \int_0^{\t^T} e^{-\rho t} L(\bar\Theta_t, l_t) dt +e^{-\rho \t^T} v(\bar\Theta_{\t^T} )\Big] 
\le v(\bar\dbx)-\frac{\e}{\rho}.
\eeaa
Since $l \in\sU_{ad}$ is arbitrary, this contradicts the dynamic programming principle (\ref{DPP_eqn}). 
%

The proof that $v$ is viscosity supersolution of (\ref{HJB_seller}) is more or less standard, again with the help of It\^o's formula (\ref{ItoForm}).
We only give a sketch here. 

Let $\bar\dbx \in \sD$ and $\f\in\ul{\sU}(\bar\dbx)$. Without loss of generality we assume that $0=v(\bar\dbx)-\f(\bar\dbx)$ is a global
minimum. That is, $v(\dbx)-\f(\dbx)\ge 0$ for all $\dbx\in \sD$. For any $h>0$ and
$l\in\sU_{ad}$,  we apply DPP (\ref{DPP_eqn}) to get
\bea
\label{super}
0&\ge& \hE\Big[ \int_0^{h} e^{-\rho t} L(\Th_t, l_t) dt + e^{-\rho {h}} v(\Th_{h}) \Big]-v(\dbx)\nonumber\\
&\ge&  \hE\Big[ \int_0^{h} e^{-\rho t} L(\Th_t, l_t) dt + e^{-\rho {h}} \f(\Th_{h}) \Big]-\f(\dbx).
\eea
Applying It\^o's formula to $e^{-\rho t} \varphi(\Th_t )$ from $0$ to $h$ we have 
\bea
\label{DCT}
 0  \geq \dbE \Big[\int_0^{h} e^{-\rho t} 
\left(L(\Th_t, l_t) +\sJ^{l} \varphi(\Th_{t})  -\rho \varphi(\Th_{t})\right) dt \Big].
\eea

Dividing both sides by $h$ and sending $h$ to $0$, we obtain
$
\rho \varphi(x,q, \hP_{\xi})
\geq \sJ^{l} \varphi(x,q, \hP_{\xi}) 
+ L(x,q, \hP_{\xi}, l) .
$
By taking supremum over $l \in \sU_{ad}$ on both sides, we conclude
\begin{equation*}
\rho \varphi(x,q, \hP_{\xi})
\geq \sup\limits_{l \in \sU_{ad}} 
[\sJ^{l} \varphi(x,q, \hP_{\xi}) 
+ L(x,q, \hP_{\xi}, l)] .
\end{equation*}
The proof is now complete.
\qed


\section{A Final Remark}
{

It is worth mentioning that the result of this paper can be connected to that of  \cite{MWZ} in the following way. 
First, note that the value function $v(x,q,\hP_\xi)$ in (\ref{value}) is the discounted lifelong expected utility of a ``representative seller",  as the
limiting case of  a Bertrand-type of game for a large number of sellers. We shall thus consider the value function of the control problem for the  representative seller as the ``equilibrium" discounted expected utility for every seller. 
Moreover, as one can see in Proposition \ref{Prop_ValueFunc}, the value function $v(x,q,\hP_\xi)$ is uniformly Lipschitz continuous, non-decreasing in $x$, and decreasing in $q$. 
Also, by Assumption \ref{assump2}-(vii), the value function is convex in $q$. 
Consequently, we see that the value function $v(x,q,\hP_\xi)$ resembles the 
%
{\it expected utility function} $U(x,q)$ in \cite{MWZ} which was defined by the following properties:

(i) the mapping $x\mapsto U(x,q)$ is non-decreasing, and $\frac{\partial U(x,q)}{\partial q} <0$, $\frac{\partial^2 U(x,q)}{\partial q^2} >0$;

(ii) the mapping $(x,q)\mapsto U(x,q)$ is uniformly Lipschitz continuous.
%

In particular, we may  identify the two functions by setting $U(x,q)= v(x,q, \hP_{\xi})|_{\xi\equiv q}$, which amounts to saying  that the 
equilibrium density function of a limit order book is fully described by the value function of a control problem of the representing seller's Bertrand-type game. This would enhance the notion of ``endogenous dynamic equilibrium limit order book model" of \cite{MWZ}  in a rather significant way. 

%
%
}


\begin{thebibliography}{99}

\bibitem{Alfonsi1}
Alfonsi, A.,  Fruth, A. and Schied, A., 
\textit{Constrained portfolio liquidation in a limit order book model}.
Banach Center Publ., {\bf 83.1} (2008), 9-25.

\bibitem{Alfonsi2}
Alfonsi, A.,  Fruth, A. and Schied, A., 
\textit{Optimal execution strategies in limit order books with general shape functions}.
Quant. Finance, 
{\bf 10}(2)  (2010), 143-157.

\bibitem{Alfonsi3}
Alfonsi, A. and Schied, A.,
\textit{Optimal Trade Execution and Absence of Price Manipulations in Limit Order Book Models}.
SIAM J. Financial Math., {\bf 1}(1) (2010), 490-522.


 \bibitem{Alfonsi4}
Alfonsi, A., Schied, A. and Slynko, A., 
\textit{Order book resilience, price manipulation, and the positive portfolio problem}.
SIAM J. Financial Math., {\bf 3}(1) (2012), 511-533.


\bibitem{Avellaneda}
Avellaneda, M. and Stoikov, S.,  
\textit{High-frequency trading in a limit order book}.
Quant. Finance, {\bf 8}(3) (2008), 217-224.

\bibitem{BHK}
Basna, R., Hilbert, A., and Kolokoltsov, V.N., 
{\it An Approximate Nash Equilibrium for Pure Jump Markov Games of Mean-field-type on Continuous State Space}.  
Stochastics {\bf 89} (2017), no. 6-7, 967-993. 

\bibitem{controlled intensity}
Bayraktar, E. and Ludkovski, M., 
\textit{Liquidation in limit order books with controlled intensity}.
Math. Finance, {\bf 24}(4) (2014), pp 627-650.

\bibitem{Bertrand}
Bertrand, J., 
\textit{Th\'eorie math\'ematique de la richesse sociale}.
 Journal des Savants, {\bf 67} (1883), 499-508.

\bibitem{Billingsley}
Billingsley, P., 
\textit{Convergence of Probability Measures}.
Second edition.
John Wiley (1999).

 
\bibitem{BLPR}
 Buckdahn, R.,  Li, J., Peng, S. and Rainer, C., 
 {\it Mean-field stochastic differential equations and associated PDEs}. 
 Ann. Probab. {\bf 45}(2) (2017),  824-878.
 

\bibitem{MFG note}
Cardaliaguet, P.,
\textit{Notes on Mean Field Games} (from P-L Lions' lectures at Coll\`ege de France).
 https://www.ceremade.dauphine.fr/ cardalia/ (2013).

\bibitem{Carmona1}
Carmona, R. and Lacker, D., 
\textit{A probabilistic weak formulation of mean field games and applications}.
Ann. Appl. Probab.
{\bf 25}(3), (2015), 1189-1231.

\bibitem{Carmona2}
Carmona, R.,   Delarue, F. and Lachapelle, A., 
\textit{Control of McKean-Vlasov dynamics versus mean field games}.
Math. Financ. Econ.
{\bf 7}(2), (2013),  131-166.

\bibitem{CS15}
Chan, P. and Sircar, R., 
{\it Bertrand and Cournot mean field games}. 
Appl. Math. Optim. 71 (2015), no. 3, 533-569.

\bibitem{Cont}
Cont, R.,  Stoikov, S. and Talreja, R.,
 \textit{A stochastic model for order book dynamics}.
 Oper. Res. {\bf 58}(3) (2010), 549-563.

\bibitem{Cournot}
Cournot, A., 
\textit{Recherches sur les Principes Mathematiques de la Theorie des Richesses}. 
Hachette, Paris, (1838). 
English translation by N. T. Bacon published in Economic Classics, Macmillan, 1897, and reprinted in 1960 by Augustus M. Kelly.

\bibitem{DI91}
Dupuis, P. and Ishii, H., 
{\it On Lipschitz continuity of the solution mapping to the Skorokhod problem, with applications}.
 Stochastics and Stochastic Reports, {\bf 35} (1991), 31-62.


\bibitem{Friedman}
Friedman, J.,
\textit{Oligopoly Theory}. 
Cambridge University Press, 1983.

\bibitem{Gatheral}
Gatheral, J., Schied, A. and Slynko, A.,
\textit{Exponential resilience and decay of market impact}.
Econophysics of Order-driven Markets, Springer (2011), pp 225-236.

\bibitem{Nadtochiy1}
 Gayduk, R. and  Nadtochiy, S.,
\textit{Liquidity Effects of Trading Frequency}.
Math. Finance 28 (3), (2018), 839-876

\bibitem{Charles}
 Gu\'eant, O., Lehalle, C.-A. and Fernandez-Tapia, J.,
\textit{Optimal portfolio liquidation with limit orders}.
SIAM J. Financial Math. 13(1), (2012), 740-764. 

\bibitem{MFG application}
Gu\'eant, O., Lasry, J.-M. and Lions,  P.-L.,
\textit{Mean field games and applications}.
 Paris-Princeton Lectures on Mathematical Finance 2010, Lecture Notes in Mathematics, vol. 2003, Springer Berlin / Heidelberg (2011), pp. 205-266.

\bibitem{MFjump}
Hao, T. and Li, J.,
\textit{Mean-field SDEs with jumps and nonlocal integral-PDEs}.
NoDEA Nonlinear Differential Equations Appl. {\bf 23}(2) (2016).

\bibitem{HHS}
Harris, C., Howison, S., and Sircar, R.,
\textit{Games with Exhaustible Resources}.
SIAM J. Appl. Math. 70(7), (2010), pages 2556-2581.



\bibitem{Ikeda}
Ikeda, N. and Watanabe, W.,
\textit{Stochastic Differential Equations and Diffusion Processes}.
North-Holland Pub. Co. (1981).

\bibitem{Kingman}
 Kingman, J. F. C.,
\textit{Poisson processes}.
Oxford University Press (1993).

\bibitem{4L}
 Lachapelle, A., Lasry, J.-M., Lehalle, C.-A. and Lions,  P.-L.,
\textit{Efficiency of the price formation process in presence of high frequency participants: a mean field game analysis}.
Math. Financ. Econ. Vol.10, Issue 3 (2016), 223-262. 




\bibitem{LS11}
Ledvina, A. and Sircar, R., 
{\it Dynamic Bertrand Oligopoly}.
App. Math. Optim.  Vol. 63, Issue 1 (2011), 11-44.


\bibitem{LS12}
Ledvina, A. and Sircar, R., 
{\it Oligopoly games under asymmetric costs and an application to energy production}. 
Math. Financ. Econ. 6 (2012), no. 4, 261-293. 

\bibitem{MFG}
Lasry, J.-M. and Lions, P.L.,
\textit{Mean field games}. 
Jpn. J. Math. 
2 (2007), 229-260. 

\bibitem{Lokka}
Lokka, A.,
\textit{Optimal liquidation in a limit order book for a risk averse investor}.
Math. Finance, Vol. 24, Issue 4 (2014), pp. 696-727. 


\bibitem{Ma93}
Ma, J., 
\textit{Discontinuous reflection and a class of singular stochastic control problems for diffusions}.
Stochastics and Stochastic Reports, Vol.44, No.3-4 (1993), 225-252. 




\bibitem{MWZ} 
Ma, J., Wang, X. and Zhang, J., 
\textit{Dynamic equilibrium limit order book model and optimal execution problem}. 
Math. Control Relat. Fields,
Vol.5, No.3 (2015). 

\bibitem{MYZ}
Ma, J., Yong, J., Zhao, Y.,  
{\it Four Step Scheme for General Markovian Forward-Backward SDEs}. 
J. Syst. Sci. Complex. 23, no. 3 (2010), 546-571.


\bibitem{thesis}
Noh, E.,
\textit{Equilibrium Model of Limit Order Book and Optimal Execution Problem}.
Ph.D. dissertation, Department of Mathematics, University of Southern California (2018).


\bibitem{Obizhaeva}
Obizhaeva, A. and Wang,  J.,
 \textit{Optimal trading strategy and supply/demand dynamics}.
Journal of Financial Markets, 71, No.4 (2013), 1027-1063.
 
 \bibitem{Potters}
Potters, M. and Bouchaud, J.-P.,
\textit{More statistical properties of order books and price impact}.
Phys. A, 324, No. 1-2 (2003), 133-140.

\bibitem{evolving_intensity}
 Protter, P., 
\textit{Point process differentials with evolving intensities}.
Nonlinear Stochastic Problems, Vol 104 of the series NATO ASI Series (1983), 133-140.

\bibitem{Protter}
Protter, P., 
{\it Stochastic Integrals and Stochastic Differential Equations}.
 (1990), Springer.


 \bibitem{Predoiu}
 Predoiu, S., Shaikhet, G. and Shreve, S.,
 \textit{Optimal execution of a general one-sided limit-order book}.
SIAM J. Financial Math. 
 2 (2011), 183-212. 

\bibitem{Rosu}
Rosu, I.,
\textit{A dynamic model of the limit order book}.
The Review of Financial Studies, 22 (2009), 4601-4641. 

\bibitem{Vives}
Vives, X.,
\textit{Oligopoly pricing: old ideas and new tools}.
MIT Press, Cambridge, MA (1999).







\end{thebibliography}
\end{document}